\documentclass[a4paper,12pt]{report}
\usepackage{graphicx,amsmath,amssymb,geometry}
\title{\bf \Huge Sheaves and Local Subgroupoids }
\author{\.Ilhan  \.I\c{c}en}
\date{}
\input{diagrams}
\def\io{^{-1}}

\newtheorem{example}{Example}[section]
\newtheorem{defn}[example]{Definition}
\newtheorem{prop}[example]{Proposition}
\newtheorem{thm}[example]{Theorem}

\newtheorem{cor}[example]{Corollary}
\newtheorem{lem}[example]{Lemma}
\newcommand{\Natno}{{\mathbb N}}
\newcommand{\Reno}{{\mathbb R}}

\newenvironment{pf}{{\bf Proof:}}{\hfill$\blacksquare$\mbox{}\vspace{.5cm}}

\def\leq{\leqslant}
\def\<{\langle}
\def\>{\rangle}
\begin{document}

\pagestyle{empty}

{~}

\vfill

\begin{center}

{\Huge \bf Sheaves and Local Subgroupoids } 

~

~

~

{\Large \bf \.Ilhan \.I\c{c}en}

\noindent
University of \.In\"on\"u,

\noindent
Science and Art Faculty,
\\
 \noindent
 Department of Mathematics\\
Malatya
/ Turkey.
\\

\noindent
E-Mail:iicen@inonu.edu.tr.          

~
~
~
~

{ \bf  University of Wales, 
Bangor, Mathematics Preprint 00.16}
\end{center}
\vspace{4ex}

{\small \noindent
{\bf Abstract:}

This is an introduction to the notion of local subgroupoid introduced
by the author and R. Brown. It can also serve as an introduction 
to an application of sheaf theory, and so could be useful to beginners 
in that theory. 

The main results are the  construction of the holonomy groupoid  
and the notion of $s$-sheaf for 
certain local subgroupoids $s$.
}

~

~
~

\vfill

~


%\include{pream}
%\pagenumbering{roman}
\tableofcontents
\pagestyle{headings}

%\pagenumbering{arabic}

%\addcontentsline{toc}{chapter}{\protect\numberline{ }{Introduction}}
\pagenumbering{arabic}
\chapter*{Introduction}

Local to global problems play a very important role in mathematics. The most important concept 
in this context is a sheaf on a topological space. A sheaf is a way 
of describing a class of functions, sets, groups, etc. For instance, a class 
of continuous functions on a topological space $X$ is very important in sheaf theory. 
The description tells the way in which a function $f$ defined on an open 
subset $U$ of $X$ can be restricted to functions $f|V$ on open subsets 
$V\subseteq U$ and then can be recovered by piecing together the restrictions
to the open subsets. This applies not just to functions, but also to other mathematical 
structures  defined `locally' on a space $X$, for example, see \cite{Ic-Yi2, St1,St2}.

Chapter 1  gives an exposition of some needed 
preliminaries and facts of the basic concepts on sheaves,  presheaves. 
An important notion in sheaf theory  is that of a global section of a 
sheaf. We show how such sections can be described
explicitly in term of 
{\it atlases}. This notion is important for all our later work. 
Also Chapter 1 reviews  the concept of local equivalence relations,  
which was  introduced by Grothendieck and Verdier \cite{Gr-Ve} in a series of exercises presented as open problems 
concerning the construction of a certain kind of topos 
was investigated further by Rosenthal  \cite{Ro1,Ro2} and more recently by  
Kock and Moerdijk  \cite{Ko-Mo1,Ko-Mo2}. 
A local equivalence relation is a global section of the sheaf ${\cal E}$ 
which is defined by the presheaf  
\[  E = \{ E(U), E_{UV}, X \}, \] 
where $E(U)$ is the set of all equivalence relations on the open subset $U$ of $X$
and $E_{UV}$ is the restriction map from $E(U)$ to $E(V)$, for $V\subseteq U$. 
Simple examples show that this presheaf is  in  genaral not a sheaf. 
It is
also becoming appaerant that groupoids are another important tool in
local-to-global problems. Therefore in Chapter 1 we also describe
basic concepts of groupoids.

Chapter 2 introduce the recent idea of  a local subgroupoid of 
a groupoid $G$ on a topological space $X$
as a global section of the sheaf ${\cal L}$ associated to the presheaf
\[            L_G = \{  L(U),  L_{UV}, X \}            \]
where $L(U)$ is the set  of all wide subgroupoids of $G|U$ 
and $L_{UV}$ is the restriction map from $L(U)$ to $L(V)$
for $V\subseteq U$.   The idea of  transitive connectness is important in the
theory and examples for local subgroupoids.
We also introduce notions of coherence which allows for an 
adjoint functional relationship between local and global subgroupoids
and  obtain a topological foliation from a local subgroupoid.

Chapter 3  defines  the holonomy groupoid of certain 
local subgroupoid by 
using the idea of locally topological groupoid.
 We define a strictly regular local subgroupoid $s$ and 
 prove that if $s$ is a strictly regular local subgroupoid
of the topological groupoid $G$ on $X$ and 
\[     glob(s) = H,  \  \  \  \  W = \bigcup_{x\in X} H_x,      \]
then $(H, W)$
admits
the structure of a locally topological groupoid. So we  obtain 
 a holonomy groupoid $H^s$ of the
strictly regular local subgroupoid $s$. 

Chapter 4  introduce the concept of $s$-sheaves for 
strictly regular local subgroupoids $s$. Corresponding concept for
local equivalence relation $r$ was extensively investigated by 
Rosental \cite{Ro1,Ro2} and Kock and Moerdjik \cite{Ko-Mo1,Ko-Mo2}
where they show that the $r$-sheaves form an \'etendue. 
This still leaves as an open problem that of describing the kind 
of topos formed by the category of $s$-sheaves.

{\bf Acknowledgements}: 
I would like to thanks to Prof.Ronald Brown,  for his suggestings, 
help  and encouragement in all stages of the  
preparation of this work and of my PhD thesis \cite{Ic2}.

\chapter{Sheaves, Atlases and Groupoids }
\section{Presheaves and Sheaves of Sets}
Let $X$ be a topological space and let ${\cal O}(X)$ be the set of open subsets of $X$.
The set ${\cal O}(X)$ is partially ordered by inclusion, so we can regard it
as a small category in the usual sense, i.e., the objects of ${\cal O}(X)$
are the open sets in $X$, and its morphisms are the inclusion maps. Also
we can form a new category ${\cal O}(X)^{op}$, called the opposite or dual
category  of ${\cal O}(X)$, by taking the same objects but reversing the
direction of all the morphism and the order of all compositions. In other
word , an arrow $V\rightarrow U$ in ${\cal O}(X)^{op}$  is the same thing
as an arrow $U\rightarrow V$ in ${\cal O}(X)$.

\begin{defn}
{\rm
Let $X$ be a topological space. A {\it presheaf} $F$ of sets on $X$ is given by
the following pieces of information;

(i) \ for each open set $U$ of $X$, a set $F(U)$,

(ii) \ for each inclusion of open sets $V\subseteq U$ of $X$, a restriction map
$F_{UV} :F(U)\rightarrow F(V)$ such that

1. \  $F_{UU} = id_U $    \   \  \   \  2. $F_{VW}\circ F_{UV} = F_{UW}$
whenever $W\subseteq V\subseteq U.$ }
\end{defn}

Thus, using functorial terminology we have the following
definition. Let $X$ be a topological space. A {\it presheaf} $F$
on $X$ is a  functor from the category ${\cal O}(X)^{op}$ of open
subset of $X$ and inclusions to the category $Sets$ of sets and
functions:
 \[           F : {\cal O}(X)^{op} \rightarrow Sets              \]

Then the system  $F = \{ F(U), F_{UV}, X \}$ is said to be  a
{\it presheaf} of sets on $X$.

In general, we define a presheaf with values in an arbitrary category.
For example, if the presheaf satisfies the following properties,
it is said to be a  presheaf
of $\Reno$-algebras:

(i) \ every $F(U)$ is an $\Reno$-algebra,

(ii) \ for $V\subseteq U$, $F_{UV}:F(U)\rightarrow F(V)$ is an $\Reno-algebra$
homomorphism.

That is, $F$ is a functor  from ${\cal O}(X)^{op}$ to the category $\Reno-Alg$ of  the
$\Reno$-algebra and $\Reno$-algebra homomorphisms:
$F : {\cal O}(X)^{op} \rightarrow \Reno-Alg $\cite{Ic-Yi1}.

Examples of presheaves  are abundant in mathematics.
For instance, if $A$ is an
abelian group, then there is the {\it constant functor} $F$ with
$F(U) = A$ for
all open $U$ and $F_{UV} = id_A$ for all $V\subseteq U$. This functor defines
the {\it constant presheaf}. We also have  the presheaf $F$ assigning to $U$ the
group (under pointwise addition)  $F(U)$ of all function from $U$ to $A$, where
 $F_{UV}$ is the canonical restriction, i.e., a functor
$F :{\cal O}(X)^{op} \rightarrow Grp.$ If  $A = \Reno $ we also have the
presheaf  $R_{\Reno}$ with  $F(U)$ being the group of all continuous real-valued
functions on $U$. Similarly, we have the presheaves of differentiable functions
on (open subsets of) a differentiable manifold $X$; of differential ${\it p-forms}$ on $X$;
of vector field on $X$; and so on. In algebraic topology, we have good
presheaf examples: the presheaf of singular $p$-cochains of open subsets
$U\subseteq X$; the presheaf assigning to $U$ its  $p$th singular cohomology group;
the presheaf assigning to $U$ the $p$th singular chain group of
$X$ $mod$ $X-U$;
and so on. For more examples, see \cite{Sw,Br,Te,St1,St2,Ma-Mo}.

Let us consider a presheaf  $F$ on a topological space as follows:
\[      F : {\cal O}(X)^{op} \rightarrow Sets.                       \]
Let $x$ be point in the  topological space $X$ and let   $U$, $V$ be
two open neighbourhoods
of $x$ and let $s\in F(U)$, $t\in F(V)$. So let us consider the set
\[        M = \{ (U,s) :  U \ \  \mbox{is open in}\ \  X, \ s\in F(U)  \}.    \]
We can define an equivalence relation on $M$ as follows; we say that $s$  and
$t$ have the {\it same germ at x} when there is some open set
$W\subseteq U\cap V$
with $x\in W$ and $s|W = t|W \in F(W)$.
The relation {\it has the same
germ at x} is an equivalence relation on $M$, and the equivalence class of any
one such  $s$ is called the {\it germ} of $s$ at $x$, in symbols $germ_xs$. Let
\[    {\cal F}_x = \{ (U, s)_x = germ_xs \mid  s\in F(U), x\in U  \ \mbox{open in} \ \ X \}           \]
be the set of all germs at $x$. Then, letting  $F^{(x)}$ be the restriction of the
functor $F : {\cal O}(X)^{op} \rightarrow Sets$ to open neigbourhoods
of $x$, the function $germ_x: F(U)\rightarrow {\cal F}_x$ forms a cone on  $F^{(x)}$
as on the right of the figure below ( because $germ_xs = germ_x(s|W)$ whenever
$x\in W\subseteq U$ and $s\in F(U)$) \cite{Ma-Mo}.

Also, if  $\{\tau_U : F(U)\rightarrow L\}_{x\in U} $ on the last below is
any other cone over $F^{(x)}$, the definition of {\it same germ} implies that
there is a unique function $t : {\cal F}_x\rightarrow L$,
with $t\circ germ_x = \tau$.
\[\begin{diagram}
     &             & F(U)    &      &                    \\
     & \SSW^{\tau_U} & \dArr & \SSE^{germ_x} &                    \\
     &              &  F(W)   &      &                    \\
     & \SW           &       & \SE  &                    \\
 L   &         & \lDotsto_{t}    &      &  {\cal F}_x
\end{diagram}\]
This just states in detail that the set ${\cal F}_x$
is the colimit and $germ_x$ is the colimitting cone of the functor $F$
restricted to open neighbourhoods of $x$:
\[        {\cal F}_x =  \lim_{\stackrel {\longrightarrow}{x\in U}}
F(U).               \]
This statement summarizes the definition of $germ$. The set $F_x$ of all germs
at $x$ is usually called the {\it stalk} of $P$ at $x$.
Now combine the various sets ${\cal F}_x$ of germs in the
disjoint union   ${\cal F}$  (over $x\in D$).
\[  {\cal F} = \bigcup_{x\in X} {\cal F}_x = \bigcup_{x\in X}\{(U,s)_x = germ_xs \mid
x\in U\subseteq X, \mbox{open}, s\in F(U) \}                                 \]
and define a canonical projection $p: {\cal F}\rightarrow X$ as the map sending
each $germ_xs = (U,s)_x$ to the point $x$, i.e., $p({\cal F}_x) = x$.

The set ${\cal F}$ will be provided with a topology such that $p$ becomes a local homeomorphism.
Let $U\subseteq X$ open and $s\in F(U)$. Then each $s\in F(U)$ determines a function
$\dot{s}$ by
\[   \dot{s} : U\rightarrow {\cal F},  \  \  \   \dot{s}(x) = (U,s)_x,  \  \  x\in U.  \]

We also define
\[   \dot{s}(U) = \bigcup _{x\in U} (U,s)_x.   \]
Topologise this set ${\cal F}$ by taking as a base of open sets
all the images $\dot{s}(U)\subseteq {\cal F}$, i.e., the family
\[ {\cal T } = \{ \dot{s}(U)  \mid U\subseteq X  \   \  open, s\in F(U) \}    \]
defines a topological base on ${\cal F}$.

Let $\dot{s_1}(U_1)$, $\dot{s_2}(U_2)\in {\cal T}$. If $\dot{s_1}(U_1)\cap \dot{s_2}(U_2) = \emptyset$,
then $\emptyset \in {\cal T}$, since $\dot{s}(\emptyset ) = \cup_{x\in \emptyset}(\emptyset, s)_x = \emptyset$.
Suppose that $\dot{s_1}(U_1)\cap \dot{s_2}(U_2)\not = \emptyset$. Then there is an element
$\sigma \in \dot{s_1}(U_1)\cap \dot{s_2}(U_2)$ such that $p(\sigma) = x\in U_1\cap U_2$.
This gives an open neighbourhood of $x\in U\subseteq  U_1\cap U_2$ such that
$\sigma = \dot{s_2}(x) = \dot{s_1}(x)$. For every $x\in U$,
since $\dot{s}(x) = \dot{s_1}(x)$,
$\dot{s}(U) = \dot{s_1}(U) \subseteq \dot{s_1}(U_1)\cap \dot{s_2}(U_2)$. So
$\dot{s_i}(U)$  lies $\dot{s_1}(U_1)\cap \dot{s_2}(U_2)$, for $i = 1,2$ and $\sigma $
is an interior point of $\dot{s_1}(U_1)\cap \dot{s_2}(U_2)$
\cite{Br3}.

Hence ${\cal F}$ is a topological space with the above topology.
This topology is called the {\it sheaf topology} on ${\cal F}$.

Now we have to show that $p$ is a local homeomorphism with this
topology, i.e., for each $\sigma = (U,s)_x\in {\cal F}$, $x\in X$,
there are open sets $U, W$ with $\sigma \in W\subseteq {\cal F}$
and
 $p(\sigma ) = x\in U\subseteq X$
such that $p|W ; W\rightarrow U$ is a homeomorphism, whereas for
$\sigma = (U,s)_x\in {\cal F}, \  p(\sigma) = p((U,s))_x = x$.
Let $\dot{s} : U\rightarrow {\cal F},  \  \  \
 \dot{s}(x) = (U,s)_x = \sigma \in {\cal F}_x$,  for $x\in U$.

Let $W = \dot{s}(U)$ and $p|U = p'$.

Firstly we will show that  $p'$ is bijective. In fact, for
$\sigma_1,\sigma_2\in \dot{s}(U) = W$, there are two elements $x_1, x_2\in U$
such that $\sigma_1 = \dot{s}(x_1)$ and  $\sigma_2 = \dot{s}(x_2)$.
If $p'(\sigma_1) = p'(\sigma_2)$, then $p'(\sigma_1) = p'(\dot{s}(x_1)) =
p'(\sigma_2) = p'(\dot{s}(x_2)) = x_1 = x_2 $. This implies $\dot{s}(x_1) =
\dot{s}(x_2)$, \ i.e., $\sigma_1 = \sigma_2$.

The map $p'$ is continuous. Choose any point $\sigma\in W = \dot{s}(U)$
such that $p'(\sigma) = x\in U$ . Then there is an open neighbourhood
$x\in V\subseteq U$ such that $\dot{s}(V)\subseteq W = \dot{s}(U)$
is an open neighbourhood of $\sigma$ and $p'(\dot{s}(V)) = V\subseteq U$.
So $p'$ is continuous.

Now we shall show that $p^{-1} = (p|W)^{-1} = \dot{s} :
U\rightarrow W = \dot{s}(U)$ is continuous. If an arbitrary
element $x\in U$, $\dot{s}(x) = \sigma \in W$, $W'\subseteq W$ is
an open neighbourhood of $\sigma $, then $(p|W)(W') \subseteq U$
is an open neighbourhood of $x$ in $U$ and  $\dot{s}(p|W) = W'$.
Hence $\dot{s}$ continuous.

These facts lead us to the basic definition of a sheaf on a topological space
$X$.
\begin{defn}{\rm
A {\it sheaf} on a topological space $X$ is a pair $({\cal F}, p)$, where

(i) \  ${\cal F}$ is a topological space (not Hausdorff in general,
see, \cite{Br3}).

(ii) \ $p : {\cal F}\rightarrow X$ is a local homeomorphism.}
\end{defn}
Then we state following theorem.
\begin{thm}\label{113}
Every presheaf $F$ on a topological space $X$ defines a sheaf ${\cal F}$
over $X$ in the above manner.
\end{thm}

\begin{defn}{\rm
Let ${\cal F}$ be a sheaf on $X$. Let $x$ be an arbitrary point in
$X$ and let $U$ be an open neighbourhood of $x$. A {\it section}
over $U$ is a continuous map $\dot{s} :U\rightarrow {\cal F}$ such
that $p\circ \dot{s} = id_U$. We denote the set of all sections of
${\cal F}$ over $U$ by $\Gamma (U,{\cal F})$. }
\end{defn}
The set of sections $\Gamma (U,{\cal F})$ gives a presheaf as follows.
If $V\subseteq U$,
\[
\Gamma_{UV}\colon \Gamma(U,{\cal F})\longrightarrow \Gamma(V,{\cal
F}),
\ \ \  \  \ \dot{s}\mapsto \dot{s}|V
\]
is the restriction, so we get a  functor
\[   \Gamma \colon O(X)^{op}\longrightarrow Sets.   \]
Hence
\[  \Gamma = \{ \Gamma (U,{\cal F}), \Gamma_{UV}, X \}     \]
is a presheaf on X. This presheaf is called the {\em canonical presheaf}.
The set of global sections of ${\cal F}$ is given by $\Gamma (X,{\cal F})$.
The presheaf  $\Gamma$
defines a sheaf $\Gamma {\cal F}$ over $X$ by Theorem \ref{113}.
Moreover
every element $s\in {\cal F}U$ is associated with a section
$\dot{s}\in \Gamma(U,{\cal F})$. If
$x\in X$ and $\sigma \in{\cal F}_x $, then there are an open
neighbourhood  $x\in U\subseteq X$ and an $\dot{s} \in \Gamma (U,{\cal F})$ such that
\[      \sigma = (U, H)_x = \dot{s}(x)  =s_x.                             \]

We shall now list some elementary properties  of the sheaf ${\cal F}$ and the
set of sections $\Gamma (U, {\cal F})$, for an open set $U\subseteq X$:

$(i)$ \ $p$ is an open map.

$(ii)$  \ Let $\dot{s} :U\rightarrow {\cal F}$ be a map
with $p\circ \dot{s} = id_U$,
for an open set $U\subseteq X$. Then $\dot{s}\in \Gamma (U, {\cal F})$ if and only if
$\dot{s}$ is open.

 $(iii)$ \ Let $U$ be open in $X$ and $\dot{s} \in \Gamma (U, {\cal F})$. Then
$p : \dot{s}(U)\rightarrow U$ is a homeomorphism and $\dot{s} =
(p|{\dot{s}(U)})^{-1}$.

$(iv)$ \  Let $\sigma$ be an arbitrary point in ${\cal F}$. Then there exists an open
set $V\subseteq X$ and a section $\dot{s} \in \Gamma (V, {\cal F})$ with
$\sigma \in \dot{s}(V)$.

$(v)$ \  For any two sections $\dot{s_1} \in \Gamma (U_1, {\cal F})$ and
$\dot{s_2}\in \Gamma (U_2, {\cal F})$, $U_1$ and $U_2$ opens, the set $U$
of points $x\in U\subseteq U_1 \cap U_2$ such that $\dot{s_1}(x) =
\dot{s_2}(x)$
is open.

Note that if  ${\cal F}$ were Hausdorff then the set $U$ of $(v)$ would also
be closed in $U_1 \cap U_2$.

1). Let ${\cal F}$ be a sheaf on a topological space $X$
and let $\Gamma $ the presheaf of sections of
${\cal F}$.  The presheaf $\Gamma $ defines a sheaf denoted by
$\Gamma {\cal F}$.
Clearly there is a natural map
\[  {\cal F}\rightarrow  \Gamma {\cal F},   \  \
 \ germ_xs\mapsto germ_x\dot{s} \]
which is a homeomorphism and preserves algebraic structure, if it has.

2). Let $F$ be a presheaf with an algebraic structure (such as
group, ring, etc.) and ${\cal F }$ the sheaf that it generates.
For such any open set $U\subseteq X$ there is a natural map
\[    \mu_U : F(U)\rightarrow \Gamma (U, {\cal F}), \ \ \  s\mapsto \dot{s}.  \]
When is $\mu_U$  is an isomorphism for all $U$ \cite{Br}. Recalling that
\[        {\cal F}_x =  \lim_{\stackrel {\longrightarrow}{x\in U}} F(U)               \]
it follows  that an element $s\in F(U)$ is in $Ker\mu_U$ if and only if $s$ is
{\it locally trivial} (that is, for every $x\in U$ there is a neighbourhood
$V$ of $x$ such that  $s|V = 0$).

Thus $\mu_u$ is a monomorphism for all $U\subseteq X$ if and only if the
following condition holds:

${\bf F_1}$  \ \ {\it If $U = \cup U_i$, with $U_i$ open in $X$,
for $i\in I$, and $s, \ t\in F(U)$
are such that $s|{U_i} = t|{U_i}$ for all $i\in I$, then $s = t$.}

Clearly, in $F_1$ we could assume that $t = 0.$ However, the condition is phrased
so that it applies to presheaf of sets.

Similarly, let $\dot{t}\in \Gamma (U, {\cal F})$. For each $x\in U$ there
is a neighbourhood $U$ of $x$ and an element $t\in F(U)_i$ with $\mu_{U_i}(t)(x)
= \dot{t}(x)$. Since $p : {\cal F}\rightarrow X$ is a local homeomorphism,
$\mu (t)$ and $\dot{t}$ coincide in some neighbourhood $V$ of $x$.
We may
assume  that $U = V$. Now
$\mu (s_i|{U_i\cap U_j}) =\mu (s_j|{U_i\cap U_j})$
so that, if ${\bf F_1}$ holds, we obtain
$s_i|{U_i\cap U_j} = s_j|{U_i\cap U_j}$.
 If \  $\Gamma $ were a presheaf sections (of any map) then this condition
would imply that the $s_i$ are restrictions to $U_i$ of a section $s\in F(U)$. Conversely,
if there is an element $s\in F(U)$ with $s|{U_i} = s_i$ for all $i\in I$, then
$\mu (s) = t$.

We have shown that, if ${\bf F_1}$ holds, then $\mu_U$ is surjective for all $U$
and (hence an isomorphism) if and if the following condition is satisfied.

${\bf F_2}$ \ \ {\it Let  $\{U_i :i\in I\}$ be a collection of open sets in
$X$ and let $U = \cup U_i$, if $s_i\in F(U_i)$ are given such that
$s_i|{U_i\cap U_j} = s_j|{U_i\cap U_j}$   for all $i, \ j$
then there exists an element $s\in F(U)$ with $s|{U_i} = s_i$ for all $i\in I$.}

Thus, sheaves are in one to one correspondence with presheaves
satisfying ${\bf F_1}$ and ${\bf F_2}$. For this reason it is
common practice not to distinguish between sheaves and presheaves
satisfying ${\bf F_1}$ and ${\bf F_2}$. Indeed, in certain
generalisations of the theory, the Definition 1.0.2 is not
available and the other notion is used. This will not be of
concern to us.

Note that ${\bf F_1}$ and ${\bf F_2}$ are equivalent to the hypothesis that the following
diagram    $(*)$  is   an equalizer diagram
$$
\begin{diagram}
  F(U) & \rDotsto^{e} & \prod_{i}{F(U)_i}  &  \pile{\rArr^{p} \\ \rArr_{q} } & \prod_{i,j}{F(U_i\cap
  U_j).}  \  \ \ \ \ \  \  (*)
\end{diagram}
$$
So we can define a ${\it sheaf}$ \  ${\cal F}$  of sets on a topological
space $X$ as a functor $F :{\cal O}^{op}\rightarrow Sets$ such
that each open covering $U = \cup U_i$, $i\in I$, of an open set $U$ of $X$
yields an equalizer diagram $(*)$, where for $t\in F(U)$,
$e(t) = \{ t|{U_i}:i\in I\}$
and for a family $t_i\in F(U)_i$, $w(t_i) = \{t_i|{U_i\cap U_j} \}$,
$q(t_i) = \{t_j|{U_i\cap U_j} \}$ \cite{Ma-Mo}.

\begin{defn}{\rm
Let $F_1$ and $F_2$ be  presheaves
of sets  on the topological space $X$.
A presheaf morphism $h : F_1\rightarrow F_2$ is a collection of
morphism $h_U : F_1U\rightarrow F_2U$ commuting with restrictions:
That is, $h$ is a natural transformation: the diagram
$$
\begin{diagram}
               F_1(U)& \rArr^{h_U}     &F_2(U)\\
        \dArr^{F_1{_{UV}}} &                 &\dArr_{F_2{_{UV}}} \\
             F_1(V)  & \rArr_{h_V}     &F_2(V)
\end{diagram}
$$
is commutative.  }
\end{defn}

\begin{defn}{\rm
Let $({\cal F}_1, p_1)$, $({\cal F}_2, p_2)$ be sheaves over $X$.

A function $\eta : {\cal F}_1\rightarrow {\cal F}_2$ is called
{\it stalk preserving} if \  $p_2\circ \eta =  p_1$
(therefore $\eta(({\cal F}_1)_x)\subseteq {({\cal F}_2)}_x$ for all $x\in X$).

A {\it sheaf morphism } is a continuous stalk preserving function
$\eta : {\cal F}_1\rightarrow {\cal F}_2$.

A {\it sheaf isomorphism} is a stalk preserving homeomorphism
$\eta : {\cal F}_1\rightarrow {\cal F}_2$. The sheaves ${\cal F}_1, \ {\cal F}_2$
are called isomorphic if there exists a sheaf isomorphism between them}.
\end{defn}

The set $Sh(X)$ will denote the category of all sheaves ${\cal F}$ of sets on
the topological space $X$ with these morphisms as arrows;
so, by definition, $Sh(X)$ is a full subcategory of the functor category
$Sets^{{\cal O}^{op}}$ \cite{Ma-Mo}.

\begin{thm}
\label{117}
Let $({\cal F}_1, p_1)$, $({\cal F}_2, p_2)$ be sheaves over $X$, and
$\eta : {\cal F}_1\rightarrow {\cal F}_2$  be a stalk preserving map. Then
the following statements are equivalent.

i). \ $\eta : {\cal F}_1\rightarrow {\cal F}_2$ is a sheaf morphism.

ii). \  For every open $U\subseteq X$ and every section
$\dot{s} \in \Gamma (U, {\cal F}_1)$, 
$\eta \circ \dot{s}\in \Gamma (U, {\cal F}_2)$.

iii). \ For every element $\sigma \in  {\cal F}_1$ there exist an open set
$U\subseteq X$ and a section $\dot{s}\in \Gamma (U, {\cal F}_1)$
with $\sigma \in \dot{s}(U)$ and
$\eta \circ \dot{s}\in \Gamma (U, {\cal F}_2).$
\end{thm}
\begin{pf}
If $\eta$ is continuous, $U\subseteq X$ open and $\dot{s}\in \Gamma (U, {\cal F}_1)$
then $\eta \circ \dot{s}$ is also continuous. Moreover
$ p_2 \circ (\eta \circ p_1) = (p_2 \circ \eta)\circ \dot{s} =
p_1\circ \dot{s} = id_U$. Therefore $\eta \circ \dot{s}$ lies in
$\Gamma (U, {\cal F}_2).$

If $\sigma \in {\cal F}_1$, then there exists an open set
$U\subseteq X$ and $\dot{s} \in \Gamma (U, {\cal F}_1)$   with
$\sigma \in \dot{s}(U)$. If the condition of (ii) are also satisfied,
then $\eta \circ \dot{s}$ lies in $\Gamma (U, {\cal F}_2)$

If for a given $\sigma \in {\cal F}_1$, open set $U\subseteq X$ and a
$\dot{s} \in \Gamma (U, {\cal F}_1)$ with $\sigma \in \dot{f}(U)$
and $\eta \circ \dot{s} \in \Gamma (U, {\cal F}_2)$ are
chosen according to condition (iii), then $\dot{s} : U\rightarrow \dot{s}(U)$
is a homeomorphism. Therefore $\eta |{\dot{s}(U)} =
 (\eta \circ \dot{s}) \circ {\dot{s}^{-1}} : \dot{s}(U)\rightarrow {\cal F}_2$ is
continuous and therefore $\eta$ is continuous at $\sigma$.
\end{pf}

\begin{cor}
Every sheaf morphism is an open map.
\end{cor}
\begin{pf}
Let $\eta :{\cal F}_1\rightarrow {\cal F}_2$ be a sheaf morphism. Since
${\cal F}_1$ is canonically isomorphic to the sheaf $\Gamma {\cal F}_1$
defined by its canonical presheaf $\Gamma_{{\cal F}_1}$. The sets $\dot{s}(U)$
with $\dot{s}\in \Gamma (U, {\cal F}_1)$ form a basis of the topology of
${\cal F}_1$. If $\dot{s}$ lies in $\Gamma (U, {\cal F}_1)$, then
$\eta \circ \dot{s}$ lies in $\Gamma (U, {\cal F}_2)$, by theorem \ref{117} (ii)
and hence
$\eta (\dot{f}(U)) = \eta \circ \dot{f}(U)$ is open in ${\cal F}_2$.
\end{pf}

{\bf Remark:}
For every open subset $U\subseteq X$,  a sheaf morphism
$\eta : {\cal F}_1\rightarrow {\cal F}_2$ defines a
map
\[  \eta_* : \Gamma (U, {\cal F}_1)\rightarrow \Gamma (U, {\cal F}_2)  \]
by  $\eta_* (s) = \eta \circ \dot{s}$.

In this stage we can define category of global sections of sheaves.

{\bf Category of Sections of Sheaves}

Let $Sh(X)$ be the category of sheaves on a topological space $X$.
We define a category of global sections of $Sh(X)$ which is
denoted by $Sec_X$ as follows;

In $Sec_X$, an arrow $\phi \colon s_1\rightarrow s_2$ is a sheaf morphism
$\phi \colon {\cal F}_1\rightarrow {\cal F}_2$ with $s_1\phi = s_2$, i.e.,
the following diagram
$$
\begin{diagram}
             {\cal F}_1 &           & \rArr^{\phi}&  & {\cal F}_2 \\
                        &\NW_{s_1}  &  &              \NE_{s_2}& \\
                        &           &       X        &  &
\end{diagram}
$$
commutes. The set of objects, $Ob(Sec_x)$ is clearly  global sections of
$Sh(X)$.

\section{Direct and Inverse Image}

Definitions of direct and inverse image of sheaves can be found any
sheaf theory book, see \cite{Sw,Ma-Mo,St1}.

Let $X$, $Y$ be  topological spaces and
let  $f : X\rightarrow Y$ be a continuous map. Then each sheaf ${\cal F}$ on $X$
yields a sheaf $f_*{\cal F}$ on $Y$ defined, for open set $V$ in $Y$
 by $(f_*{\cal F})(V) =
{\cal F}(f^{-1}V)$; that is, $f_*{\cal F}$ is defined as the
composition functor
\[   {\cal O}(Y)^{op} \stackrel{f^{-1}}\longrightarrow {\cal O}(X)^{op}
\stackrel{\cal F}\longrightarrow Sets   \ \   \]

This sheaf $f_*{\cal F}$ is called the {\it direct image of ${\cal F}$ under $f$}.
The map $f_*$ so defined is clearly a functor
\[ f_*  : Sh(X) \rightarrow Sh(Y).  \]
Also $(fg)_* = f_* g_*$, so the definition $Sh(f) = f_*$ makes $Sh$ a functor on
the category of all  topological spaces. In particular, if $ f :X\rightarrow Y$
is a homeomorphism, $f_*$ gives an isomorphism of categories between sheaves on $X$
and on $Y$.

Now let ${\cal F}$ be a sheaf on $Y$. The {\it inverse image} $f^*{\cal F}$ of $F$ is the sheaf
on $X$ defined by
\[ f^*{\cal F} = \{ (x, \sigma) \in X\times {\cal F} : f(x) = p(\sigma)  \}    \]
where $p : {\cal F}\rightarrow Y$ is the canonical projection of sheaf,
i.e. $p$ is a local homeomorphism. A projection on $f^*{\cal F}$ is given by
\[  p^* :  f^*{\cal F}\rightarrow X,  \   \   \ (x, \sigma )\mapsto x. \]
To check that $f^*{\cal F}$ is indeed a sheaf we note that if $U\subseteq Y$ is an
open neighbourhood of $f(x)$ and $ \dot{s} : U\rightarrow {\cal F}$ is a section
of ${\cal F}$ with $\dot{s}(f(x)) = \sigma$, then the neighbourhood $(f^{-1}(U))\times
\dot{s}(U))\cap f^*{\cal F}$ if $(x, \sigma )$ in $f^*{\cal F}$ is precisely
\[     \{ (x', \dot{s}f(x'))\colon  \  \  x'\in f^{-1}(U)  \}  \]
and hence maps homeomorphically onto $f^{-1}(U)$.

Then, each continuous map $f :X\rightarrow Y$ gives a functor
$f^* :Sh(Y)\rightarrow Sh(X)$. For a sheaf ${\cal F}$ on $Y$,
the value $f^*{\cal F}\in Sh(X)$
of this functor is called the {\it inverse image of ${\cal F}$ under $f$}.

\begin{thm}
If $f : X\rightarrow Y$ is a continuous map, then the functor $f^*$,
sending each  sheaf ${\cal F}$ on $Y$ to its inverse image on $X$, is left
adjoint to the direct image functor $f_*$;
$$
\begin{diagram}
           Sh(X)& \pile{\rArr^{f_*} \\ \lArr_{f^*}} & Sh(Y)
\end{diagram}
$$
\end{thm}
\begin{pf}
See Mac Lane  and Moerdijk \cite{Ma-Mo}.
\end{pf}

\begin{defn}{\rm
Let ${\cal F}$ be a sheaf on $X$ and  $Y\subseteq X$. Then
\[    {\cal F}\mid_Y = p^{-1}(Y)                   \]
is a sheaf on $Y$ called the {\it restriction of ${\cal F}$} on $Y$.}
\end{defn}

Let $F$ be a constant presheaf on $X$. A sheaf ${\cal F}$
which is generated by $F$ is called {\it constant sheaf}.
In other word, the constant
sheaf on $X$ with  stalk $F(U) = A$ is the sheaf $X\times A$
(giving $A$ the
discrete topology).

A sheaf ${\cal F}$ on $X$ is called {\it locally constant}
if every point of $X$
has a neighbourhood $U$ such that ${\cal F}\mid_U$ is constant.
If a presheaf $F$ is a sheaf in the  functorial terminology,
then the locally constant sheaf can be defined as
follows;
A sheaf ${\cal F}$ on a topological $X$ is called {\it locally constant}
if each point $x\in X$ has a basis of open neighbourhood $N_x$ such that
whenever $U, V\in N_x$ with $V\subseteq U$, the restriction
\[             F_{UV} : F(U)\rightarrow F(V)                  \]
is a bijection.
We also say that ${\cal F}$ is {\it locally constant}
if and only if $p : {\cal F}\rightarrow X$
is a covering  \cite{Ma-Mo}.

We will  give  definition of {\it atlas} and  {\it chart} due to \cite{Ko-Mo1},
for any presheaf $F$ or,
more precisely, an {\em atlas} for a global section of the sheaf  ${\cal F}$
associated to a presheaf $F$.

\section{Chart and Atlas}

A {\it global section} of the sheaf  ${\cal F}$
associated to a presheaf $F$ on a topological space $X$ can be constructed
in  different ways. In sheaf theory, this is usually done by constructing
the sheaf space, consisting of germs of elements of $F$, at
various points, but for the present purpose, a description in term of
atlases is more appropriate \cite{Ko-Mo1}.

Suppose we are given a presheaf $F$ on a topological space $X$.
An {\it atlas} in the presheaf $F$,
or an atlas for a global section of ${\cal F}$
, consists of a family,
\[   {\cal U} = \{ (U_i, s_i) : i\in I, \  s_i\in F(U_i) \}    \]
such that

(i) \ \  $X = \displaystyle \bigcup_{i\in I} U_i$ ,i.e., the family  $\{U_i : i\in I\}$ is an open covering of $X$.

(ii) \ \ ({\it Local compatibility Condition}). For all $i, j\in I$, there
exist an open cover of $U_i\cap U_j$  by open sets $W$ for which
\[       s_i|W = s_j|W.                                \]

If ${\cal U}$ is an atlas as above,
 then  each $(U_i,s_i)$ is called its  {\it chart}.

\begin{prop}
Every global  section  $s$  of  the  sheaf  ${\cal F}$ associated
to  presheaf  $F$ is  given  by  some
atlas. Conversely,  every atlas  in  $F$   determines  a  global  section.
\end{prop}
\begin{pf}
Let
\[   {\cal U} = \{ (U_i, s_i) : i\in I, \  s_i\in F(U_i) \}    \]
be an atlas for a presheaf $F\colon {\cal O}(X)^{op}\rightarrow Sets$.

We claim that the atlas   ${\cal U}$  defines a global section of the sheaf
${\cal F}$  associated to the presheaf $F$ above.

Clearly we can define an equivalence relation on the atlas ${\cal U}$
as follows; let us fix an element $x\in X$. Let $(U_i, s_i)$,
$(U_j, s_j)$ be two elements of ${\cal U}$ such that $x\in U_i\cap U_j$.
We say that  $(U_i, s_i)$, $(U_j, s_j)$ are equivalent, written
$(U_i, s_i) \sim_x (U_j, s_j)$ iff there is a neighbourhood $W$ such that
$x\in W\subseteq U_i\cap U_j$ and $s_i|W = s_j|W$.
Let $(U_i, s_i)_x$ denote
the equivalence classes of $(U_i, s_i)$. Then we obtain stalks and their sheaf
as usual.
\[  {\cal F}_x = \{ (U_i, s_i)_x : x\in U_i, \  s_i\in F(U_i) \},
\  \  \  {\cal F} = \displaystyle \bigcup_{x\in X} {\cal F}_x    \]

So each $(U_i, s_i)_x$ determines a map
$\dot{s_i}$ by
\[   \dot{s_i} : U_i\rightarrow {\cal F}  \  \  \
\dot{s_i}(x) = (U_i,s_i)_x  \  \  x\in U \]
such that $\dot{s_i}$ is continuous in the sheaf topology of ${\cal F}$.
Since ${\cal U}$ is an atlas, then the formula
\[   \dot{s}(x) = \dot{s_i}(x)   \  \mbox {for}  \  \  \  x\in U_i.  \]
defines a map $\dot{s}$ of the topological space $X$ into  the sheaf ${\cal F}$.
For an arbitrary $U$ which is open in ${\cal F}$ we have
\[     \dot{s}^{-1}(U) = \bigcup_{i\in I}\dot{s_i}^{-1}(U).  \]
The set $\dot{s_i}^{-1}(U)$ is open in $U_i$, whence also in $X$. It follows that
$\dot{s}^{-1}(U)$ is an open set in $X$. Hence
$\dot{s} : X\rightarrow {\cal F}$ is continuous.

In addition, we have to show that $p\circ \dot{s} = id_X$, where $p\colon {\cal F}\rightarrow X$
is a local homeomorphism.
For any $x\in X$ we have an open set $U_i$ containing $x$. So
$p\circ \dot{s}(x) = p\circ \dot{s}_i(x) = p((U_i, s_i)_x) = x$. Thus $\dot{s}$
is a global section of the sheaf ${\cal F}$.

Conversely, a global section of a sheaf defines an atlas.

Let $\dot{s}$ be a global section of a sheaf ${\cal F}$ associated to the presheaf $F$.
This means that there exists a continuous map $\dot{s}\colon X\rightarrow {\cal F}$
such that $p\circ \dot{s} = id_X$, where $p\colon {\cal F} \rightarrow X$
is a local homeomorphism.
Since $\dot{s}$ is continuous, each point $x\in X$ has an open
neighbourhood  $ U_i$
such that $\dot{s}|{U_i}$ is continuous.  Let $\dot{s}|{U_i} = \dot{s_{i}}$.
Then  we have a continuous map
 $\dot{s_i}\colon U_i\rightarrow {\cal F}$ defined by
 $\dot{s_i}(x) = (U_i, s_i)_x$,
where $(U_i, s_i)_x$ is the equivalence class of $(U_i, s_i)$.
That is, each $x\in X$ has an open set $U_i$ with $x\in U_i$,
and every map $\dot{s}_i$  over $U_i$ defines $(U_i, s_i)$ .
Indeed, these  $(U_i, s_i)'s$
form  the following atlas:
\[    {\cal U} = \{ (U_i, s_i) : i\in I, s_i\in F(U_i) \}.   \]
\end{pf}
An atlas ${\cal V} = \{(V_i,t_i) : i\in I\}$  is called a {\em refinement} of  ${\cal U}$ if there
is a function $p\colon J\rightarrow I$, such that, for each
index $V_j\subseteq U_{p(j)=i} $
and $s_i|{V_j}  = t_j$ , i.e. $(V_j,t_j)$ is sub-chart of $(U_i,s_i)$, $j\in J$.
Moreover, two atlases define  the same global section $s$ if and only  if   they
have common refinement \cite{Ko-Mo1}.

Given a global section $s$, a pair $(U,t)$ , where  $t\in F(U)$,  may  be
 called a {\it chart} for s if either of the two equivalent  conditions
hold ;

i)  There exist some atlas $U$ for $s$ with $(U,t)$ as a member.

ii) For every $x\in U$ the germ of $t$ at $x$ equals $s(x)$.

\section{Local equivalence relations}
Recall that  the notion of local equivalence relation on a
topological space $X$ was introduced by  Grothendieck and Verdier
\cite{Gr-Ve}
 p.485 in series of exercises, presented essentially as open problems
for the purpose of constructing a certain kind of topos called an
\'etendue. The concept has been investigated  by Rosenthal
\cite{Ro1,Ro2} and more recently by Kock and Moerdijk
\cite{Ko-Mo1, Ko-Mo2}.

Let $X$ be a topological space and let $U$ be an open subset of $X$. Let $E(U)$ be the set of all equivalence relations
on $U$. If $V\subseteq U$ is also open,
we have a restriction morphism
\[
\begin{array}{cccl}
E_{UV}\colon& E(U)&\longrightarrow &E(V)   \\
            &    R & \mapsto &  R|V =R\cap (V\times V)
\end{array} \]
defines a  functor
\[ E\colon {\cal O}(X)^{op}\longrightarrow Sets  \]
 from the category of open subsets of $X$ and
inclusions to the category of sets and functions.
Hence $E$ becomes a presheaf on $X$. This presheaf is denoted by
\[  E = \{ E(U), E_{UV}, X \}. \]
Rosenthal showed  in his paper that this
presheaf is  in general not a sheaf \cite{Ro1}.
In fact, if $E$ is a sheaf, it must satisfy the  conditions $F_1$ and
 $F_2$
in first section.

We now give  Rosenthal's example.
\begin{example}{\rm
Let $X\subseteq {\Reno^2}$ be defined by
\[ X=\{(0,y) : 0\leq y\leq 1\}\cup \{(x,1) : 0\leq x \leq 1 \}\cup\{V_n : n\in {
\Natno }\} \]
where $ V_n = \{(1/n, y) : 0\leq y \leq 1 \}$ (See Figure 1).

Let $\{ U_i : i\in I \}$ be an open covering of $X$.
We can define an equivalence relation $R_i$ on each open set
$U_i$, for $i\in I$ as follows:
\[    x R_i  y   \Leftrightarrow   \mbox{there is a path in $U_i$ joining $x$ to $y$. }  \]
Let $E(U_i)$ denotes the set of all equivalence relations on $U_i$. This lead us to
the definition of presheaves defined by equivalence relations as mentioned earlier.
Namely,
\[ E\colon {\cal O}(X)^{op}\rightarrow Sets \]
is a functor.

Now, let us choose open sets $\{ U_1, U_2, U_3 \}$ of $X$  such that
$U = \displaystyle \bigcup_{i=1}^3 U_i$  as in the following
diagram:

\begin{center}
\begin{figure}[ht]
\begin{center}
\includegraphics[angle=-90]{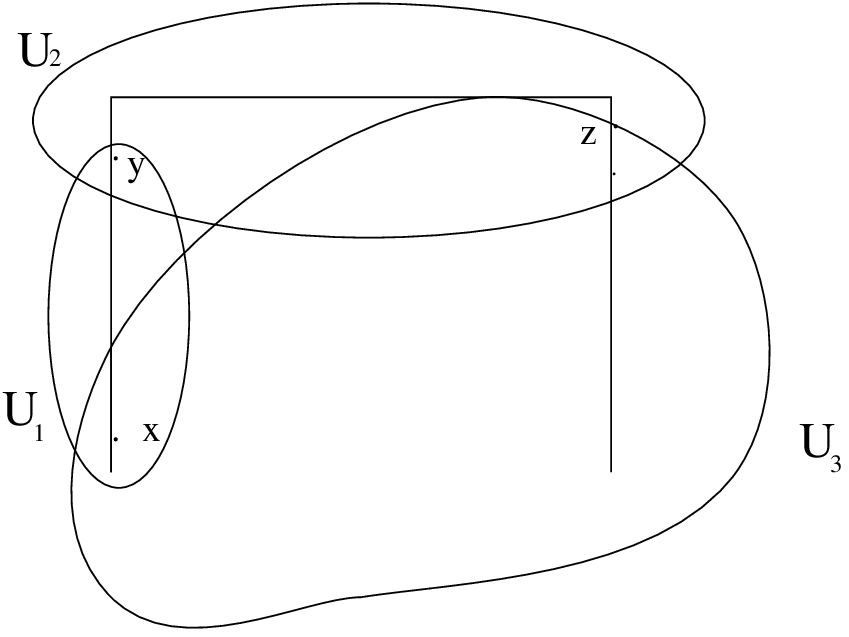}
\end{center}
\end{figure}
\end{center}
\newpage
\begin{center}
Figure 1
\end{center}

The presheaf $E$ is easily seen to satisfy the condition
${\bf F}_1$ for a sheaf. Now we have to
show that presheaf $E$ must satisfy the second condition ${\bf F}_2$, see
previous sections, i.e., there should be an equivalence relation $R\in E(U)$ such that
\[                R|{U_i} = R_i,   \  \mbox{for} \ \ \ \  i = 1, 2, 3.  \]
However it does not satisfy $F_2$. Suppose we are given such an
equivalence relation $R$ on $U$. Let us choose arbitrary points
$x, y, z $ as in the diagram. Since $R|{U_1} = R_1$, then there is
a path in $U$ joining $x$ to $y$. Also $R|{U_2} = R_2$, then there
is a path in $U$ joining $y$ to $z$. But $x, z\in U_3$, so $R_3$
should have to satisfy that there is a path in $U_3$ joining  $x$
to $z$. This is a contraction, because there is no such a path
joining them in $U_i$. So $R|{U_3} \neq R_3$. This shows that $E$
is not a sheaf on $X$. }
\end{example}

In first section, we have showed that  a sheaf can be
constructed for each presheaf. Hence, let
${\cal E}\rightarrow X$  denote the sheaf  corresponding to  $E$.
Let $U\subseteq X$ be open and let $r\colon U\rightarrow {\cal E}$ be a section
of ${\cal E}$ over $U$.
The set of such sections is denoted by $\Gamma(U,{\cal E}) $.
The set $\Gamma (X, {\cal E}) $ is  the set of global sections
of ${\cal E}$.
Also $\Gamma $ is a presheaf.

\begin{defn}
{\rm
A global section $r$ of the sheaf ${\cal E}$ is called a {\it local equivalence relation }
on $X$.
}
\end{defn}

A local equivalence relation $r$ may be given by the following
local data which is called an {\it atlas}: For
an open cover $\{U_i:  i\in I\}$ of $X$, if $ R_i\in E(U_i)$,
$R_j\in E(U_j)$
such that
 $z\in U_i\cap U_j$, there is an open neighbourhood $W$ of $z$ with
$W\subseteq U_i\cap U_j$ and $R_i|W = R_j\mid _W$. This condition
is called the { \it local compatibility condition}.
Conversely, by our earlier discussion, every local equivalence
relation  $r$ is defined by an atlas
which is denoted by
${\cal U}_r= \{(U_i, R_i): i\in I, R_i\in E(U_i)\}.$

In his paper, Rosenthal gives many examples of local
equivalence relations.
His  main example comes from foliation.

Let $X$ be a $C^{\infty }$- manifold of dimension $n$. A {\it foliation } of
codimension $q$ (dimension $p$, where $p+q = n$) is given by the following;

(1) an open cover $ \{U_i:i\in I \}$  of $X$.

(2) submersions $f_i \colon U_i\rightarrow \Reno^q $ for $i\in I$. That is,
 for all $i\in I$, $f_i$ is smooth and for all $x\in U_i$,
 there is an open neighbourhood $V$ of $f_i(x)$
and a smooth section $g\colon V\rightarrow U_i $ with
$g(f_i(x)) = x$.

(3) for $i,j\in I $ and $x\in U_i\cap U_j $, there is a diffeomorphism $\gamma_{j,i}$
from a neighbourhood of $f_i(x)$ such that locally $f_i = \gamma_{j,i} f_i$  and
if $x \in U_i \cap U_j \cap U_k$, then locally
\[ \gamma_{k,i} = \gamma_{k,i} \gamma_{j,i}  \]
\begin{example} {\rm
A foliation on a manifold makes $X$ look locally like $n$-space divided into
parallel $p$-planes. The simplest example of a foliation is given by the submersion
$$
\begin{diagram}
 \Pi_2\colon { \Reno^p\times \Reno^{n-p}\longrightarrow  \Reno^{n-p}  }
\end{diagram}
$$
and every foliation locally looks like this. On each $U_i$, let $R_i$ be defined
by $(x,y)\in R_i$ iff $f_i(x) = f_i(y)$. Condition (3) above is exactly
the local compatibility of the equivalence relations.
Let $U_i$ and $U_j$ be  open sets in $X$ and let $R_i$ and $R_j$ be
corresponding equivalence relations on $U_i$ and $U_j$, respectively.
For any $x\in U_i\cap U_j$, there is a open set $U_k$ such that
$x\in U_k\subseteq U_i\cap U_j$.
Choose two points $x,y\in U_k$ such that $(x,y)\in R_i$ , so $f_i(x)=f_i(y)$.
By the definition of foliation, $\gamma_{j,i}(f_i(x))= f_j(x)$ and
$\gamma_{j,i}(f_i(y)=f_j(y)$. Since $\gamma_{j,i}$ is a diffeomorphism,
$f_i(x) = f_j(y)$. Then $(x,y)\in R_j$ and so $R_i|{U_k} = R_j|{U_k}$.
Hence we get a local equivalence relation
$r$ on $X$. For more information about and examples of foliations,
see  \cite{Ca-Co,La1,La2,Mo2,Ph,Wi}.}
\end{example}

Suppose we are given topological spaces $X$ and $Y$, an open
 cover $\{ U_i\colon i\in I\}$
of $X$ and a family of mapping  $\{f_i\colon i\in I \}$ , where $f_i\colon U_i\rightarrow Y$, and
\[          f_i|{U_i\cap U_j} = f_j|{U_i\cap U_j}             \]
for every $i,j\in I$. If the mappings are compatible, then the formula
$f(x)=f_i(x)$ for $x\in U_i$ defines a function $f$ from the space $X$ into the
space $Y$ such that $f|{U_i}=f_i$ \cite{Bo2}.

\begin{example}{\rm  So the above example, clearly,
could be generalised to a space $X$ and a cover $\{U_i\}$
and  continuous $f_i\colon X\rightarrow Y$ for some space $Y$ with the
necessary local compatibility.
We define  an equivalence relation $R_i$ on $U_i$ by  $xRy$ if and only if
 $f_i(x)=f_j(x)$. Let $R_i$ and $R_j$ be the equivalence relation on $U_i$ and
$U_j$, respectively, and let $x\in U_k\subseteq U_i \cap U_j$ for $x\in U_i\cap U_j$.
By local compatibility condition, $R_i|{U_k} = R_j|{U_k}$.}
\end{example}

We now give another example of a local equivalence relation using the notion
of sober space.

A topological space $X$ is said to be  ${\it sober }$ iff for any open subset
$U\in O(X)$ such that

(i) \ $U\neq X$

(ii) \ for all open $W,V$ in $X$, if $W\cap V \subseteq U$,

then $W\subseteq U$ or $V\subseteq U$,
then there is a unique point
$x\in X$ with  $U = X-\{x\}$.

Again  Rosenthal's paper gives an example of a topological space
which is not sober.

This definition is often phrased in terms of closed sets. A closed subset
$F\subseteq X$ is called {\it irreducible} if it can not be written as the union of two smaller
closed subsets; that is, whenever $F_1$ and $F_2$ are closed sets with
$F = F_1\cup F_2$ , then $F_1 = F$ or $F_2 = F$. Clearly, if $y$ is a point of $X$
then $\overline {\{ y\}}$ is an irreducible closed set. Thus $X$  is
{\it sober}
iff every non empty irreducible closed set is the closure of a unique point.

For observe that, for any open set $U\subseteq X$ and its closed complement $F = X-U$
, the set $U$ is proper prime, as in (i) and (ii), iff
$F$ is non empty and irreducible.

Any Hausdorff space $X$ is sober; any sober space is $T_0$.  Sobriety is a weaker separation axiom than
$T_1$. An important example of sober spaces in algebraic geometry is the
Zariski spectrum of a commutative ring \cite{Ma-Mo}.
\begin{example}{\rm
If we take a space $X$, and an open set $U\subseteq X$, we have the equivalence relation
$x\sim y$ iff $\overline {\{x\}}^U = \overline {\{y\}}^U $, where we take
closure relative to $U$. This  equivalence relation defines a local equivalence relation on $X$.}
\end{example}

\section{The Category of Groupoids}

The rest of this section is to give some knowledge of groupoid theory,
as contained in \cite{Br3,Hi,Mac,Mo2}.

An interesting algebraic theory of groupoids exists, and was begun
by Brandt and by Baer in the 1920's,  well before Ehresmann
introduced the concept of groupoid into differential geometry and
topology in the 1950's.

We give some references for the algebraic theory of groupoids and their application to problems in
group theory see \cite{Hi}, and  for general topological groupoids, see
\cite{Br-Hr,Mac}.

The discussion is carried only as far as is needed to supply a
convenient topological framework for the discussion throughout
rest of the study.

We begin by reviewing
the basic definitions and fixing the notations.
Recall that simply  a {\bf groupoid} is a small category in which each
morphism is an isomorphism. Here we shall give explicit definition of the
groupoid and properties.

\begin{defn}{\rm
A {\it groupoid} is a category $G = ( G, X, \alpha, \beta, m, i )$
given by a set $G$ of arrows, a set $X$ of objects and four structure maps;
$$
\begin{diagram}
G\times_XG  & \rArr^{m}  & G & \pile{\rArr^{\alpha ,\beta}\\ \rArr \\
\lArr_{i} } & G
\end{diagram}
$$
The maps $\alpha$ and $\beta$ give for each arrow $g\in G$ its {\it source}
  (domain) $\alpha (g)$ and its {\it target} (codomain) $\beta (g)$.
The map $m$ is defined for any pair of arrows $f, g$ with $\alpha
(f) = \beta (g)$, and assigns to this pair the composition $(f,
g)$ also denoted $f\circ g$. Finally, the map $i$, called
inclusion map, assigns to each object $x\in X$ the identity arrow
at $x$, denoted $i(x)$ (or $id_x$ or $1_x$). These maps must
satisfy the well-known identities.

$\alpha (i(x)) = x = \beta (i(x))$,   \  \  \  \ $(f\circ g)\circ h = f\circ (
g\circ h),$

$\alpha (f\circ g) = \alpha (g)$,   \  \  \  $f\circ i(\alpha (f)) = f,$

$\beta (f\circ g) = \beta (f)$,  \  \  \  $i (\alpha (f)) \circ = f,$

\noindent
and for each $f\in G$, $f\colon x\rightarrow y$, there exist an arrow
$g\colon y\rightarrow x$, $g\in G$, so that $f\circ g = i(y)$ and
$g\circ f = i(x)$.}
\end{defn}

\
\
\

Intuitively, one thinks of a groupoid $G$ as a disjoint union of
the set of arrows $G(x,y)=\{f\in G\mid f:x\rightarrow y\}$
parameterized by $x,y\in X$. Namely,
\[    G = \bigcup_{x, y\in X} G(x, y).      \]
\begin{defn}{\rm
A  category  $G'=( G', X', \alpha ', \beta ', m', i')$ is a subgroupoid
of a groupoid $G = (G, X,\alpha, \beta, m, i)$  provided $G'$ is a subset of arrows $G$,
and $X'$ is a subset of objects $X$ and its four
structure maps are restrictions and $G'$ is a groupoid.}

\end{defn}

A subgroupoid  $G'$ is called {\it full}({\it wide}) subgroupoid if $G'$
is a full(wide) subcategory. That is , the subgroupoid  $G'$
is full(wide) subgroupoid if $G = G'$ $( X = X' )$, respectively.
A groupoid $G = ( G, X, \alpha, \beta, m, i)$ is said to be
{\it transitive}
if $G(x,y)=\{f\in G\mid f\colon x\rightarrow y , \forall x, y\in X \}$ is
non-empty and {\it totally transitive} if the set
$G(x,y) = \{f\in G\mid f\colon x\rightarrow y  \}$ has a single element and
a groupoid is said to be {\it locally transitive} if $X$ has a basis of
open sets $U$ such that the restriction of $G$ to $U$ is transitive,
similarly for {\it simply transitive} and {\it locally simple transitive}, so on.

Let $G = ( G, X, \alpha, \beta, m, i )$ be a groupoid and let $a\in X$.
Let $T_a$ be the full subgroupoid of $G$ on all objects $y\in X$ such that
$G(a, y)$ is non-empty. Then, if $x, y\in Ob(T_a)$, $G(x, y)$ is non-empty,
since it contains the composite $gf$ for some $f\in G(x, a)$ and some
$g\in G(a, y)$. Thus $T_a$ is transitive and is clearly the maximal
transitive subgroupoid of $G$ with $a$, as one of its object $T_a$ is,
therefore called the {\it component} of $G$ containing (or determined by )
$a$ \cite{Br3}.
\subsection{Examples}

\begin{example} \label{143}{\em
Any set $X$ may be regarded as a groupoid on itself with $\alpha = \beta = id_X$
and every element an identity. Groupoids in which every element is an identity
have been given  a variety of names, we will call them {\bf null} groupoids.
It has only identities $1_x$, one for each $x\in X$, these may be composed with
themselves so that $1_x . 1_x = 1_x$ and there are no other composition.}
\end{example}
\begin{example}\label{145}{\em
Let $X$ be a set and then there is a groupoid with object set $X$
and set of arrows the product set $X\times X$ so that an arrow
$x\rightarrow y$ is simply the ordered pair $(y, x)$. The composition is
then given by
\[     (z, y) (y, x) = (z, x).   \]
This groupoid looks rather simple, banal and unworthy of
consideration. Surprising, though, it plays a key role in the
theory and application. One reason is that if $G$ is a subgroupoid
of $X\times X$ and $G$ has the same object set $X$ as $X\times X$,
then $G$ is essentially an equivalence relation on $X$. That is,
for all $x\in X$;
 $(x, x)\in G$; if $(x, y)\in G$ then $(y, x)\in G$; and if
$(z, y), (y, x)\in G$ then $(z, x)\in G$.
Now equivalence relation is a groupoid with the composition above, is
important in mathematics
and science because they formalize the idea of classification -two elements
have the same classification if and only if they are equivalent. In
mathematical terms, we say that equivalence relations formalize the idea of
quotienting. Thus, it is an important aspect of their applications
that groupoids generalise both groups and equivalence relations \cite{Br2}.}
\end{example}

\begin{example}\label{group}{\em
Let $X$ be a set and $K$ is a group. We give  $ X \times X\times K$
the
structure of a groupoid on $X$ in the following way,
 $\alpha$ is the projection onto
the second factor of  $ X\times X\times K$
and $\beta$ is the projection onto the first
factor:
\[     \alpha (x, y, g) = y    \   \   \  \  \beta (x, y, g) =x   \]
the inclusion map is $x\mapsto 1_x = (x, x, 1)$ and the composition is
\[  (z, y, h)(y, x, g ) = (z, x, hg )    \]
defined iff $ y = y'$. The inverse of $(y, x, g)$ is $( x, y, g^{-1})$.
         }
\end{example}

\begin{example}{\em  \cite{Ic-Yi2}
A {\em pointed space} is a pair $(X,x)$ where $x\in X$ and $X$ is a topological
space. A {\em pointed map}  $(X,x)\rightarrow (X,y)$ on $X$ is determined by  two
pointed spaces and a map $\psi\colon X\rightarrow X $ such that $\psi x = y$.
We get a category $Top_*$ of pointed spaces (or spaces with base point). We shall define
an equivalence relation in this category: A map $\psi\colon (X,x)\rightarrow (X,y)$
called  a {\em homotopy equivalence} if there exists a  map
$\psi^{-1}\colon (X,y)\rightarrow (X,x)$ such that
$\psi^{-1}\psi\simeq 1_{(X,x)}$ and $\psi\psi^{-1}\simeq 1_{(X,y)}$, where
$\simeq$ is the homotopy relation (rel base points).
Clearly this relation is an
equivalence relation on $Top_*$. Let $[(X,x);(X,y)]$ denote the
set of all homotopy  classes of homotopy equivalence
maps $(X,x)\rightarrow (X,y)$.

We consider    a    groupoid    over   $X$,    called  ${\cal E}(X)$,    such    that
${\cal E}(X)(x,y)=[(X,x),(X,y)]$ is the set of pointed  homotopy
classes
of homotopy equivalence maps  $(X,x)\rightarrow (X,y)$. So the set
\[              {\cal E}(X) = \bigcup_{x,y\in X}[(X,x);(X,y)]           \]
is a groupoid under the composition: ${\cal E}(X)\oplus {\cal E}(X)\rightarrow {\cal E}(X)$,
$([\psi],[\psi'])\mapsto [\psi][\psi'] = [\psi \psi']$, where
${\cal E}(X)\oplus {\cal E}(X) = \{([\psi],[\psi'])\colon \beta [\psi'] =\alpha [\psi] \}$.
For any element $[\psi]\in {\cal E}(X)(x,y) = [(X,x);(X,y)]$, the  source and target maps are
$\alpha [\psi] = x$ and $\beta [\psi] = y$ respectively and $\varepsilon \colon X\rightarrow {\cal E}(X)$
, $x\mapsto [1_x]$.}
\end{example}
\begin{example}\label{147}{\em
An interesting  groupoid is the {\it fundamental groupoid}  $\pi_1
(X)$ of a topological space $X$. An object of $\pi (X)$  is a
points $x$ of $X$, and arrow $a \colon x\rightarrow y$ of $\pi
(X)$ is a homotopy class $[a]$ of paths $a\colon [0,1]\rightarrow
X$ from $x = a(0)$ to $y = a(1)$. Such a path  $a$ is a continuous
map $I = [0, 1] \rightarrow X$ with $a(0) = x$, $a(1) =x'$ while
two paths $a, b$ with the same end points $x$ and $x'$ are
homotopic, when there is a continuous map  $F\colon I\times
I\rightarrow X$ with $F(t, 0) = b(t)$, $F(t, 1) = a(t)$, and $F(0,
s) = x$, $F(1, s) = x'$, for all $s$ and $t$ in $I$. The
composition  of paths $a\colon x\rightarrow x'$ and $b\colon
x'\rightarrow x''$ is the path $c$ which is {\it a followed by g}
given explicitly by
\[
h(t)  =  \left\{\begin{array}{cl}
            a(2t)  ,&  \hspace{1cm}          0\leq t\leq 1/2 \\
            b(2t-1),&  \hspace{1cm}       1/2\leq t\leq  1
\end{array}   \right.   \]
Composition applies also to homotopy classes, and makes $\pi (X)$ a groupoid.
\cite{Br3,Mo2,Ma}.}
\end{example}

Also the category of finite set and bijections, and an individual group
can be given as examples of groupoids \cite{Br4}.

\begin{defn}{\rm
Let $G_1 = ( G_1, X_1, \alpha_1, \beta_1, m_1, i_1 )$  and
 $G_2 = ( G_2, X_2, \alpha_2, \beta_2, m_2, i_2 )$  be two groupoids.
A groupoid morphism
\[ \phi \colon G_1 \rightarrow G_2  \]
is a functor, i.e., a groupoid morphism between two groupoids
$G_1$  and
$G_2$ is given by two morphisms
an objects $Ob(\phi)$ and arrows $\phi$ (both will be denoted $\phi$):
\[  Ob(\phi ) = \phi \colon X_1\rightarrow X_2,    \  \   \
\phi \colon G_1\rightarrow G_2                  \]
with respect the structure maps of the groupoid.

$\phi (i(x)) = i(\phi (x))$, for  each $x\in X$,

$(\alpha_i(g)) = \alpha_i (\phi (g)),$ \  \  $\phi (\beta_i (g)) = \beta_i (\phi
(g)),$
$i = 0, 1$, for  each $g\in X_1$.

$\phi (f\circ g) = \phi (f)\circ \phi (g)$ all $ f,g \in  G_1$ \ \  \mbox{with}
\ \
 $\alpha_i(f) = \beta_i (g)$.}
\end{defn}.
\begin{defn}{\rm
Let $G_1 = ( G_1, X, \alpha_1, \beta_1, m_1, i_1 )$   and
 $G_2 = ( G_2, X, \alpha_2, \beta_2, m_2, i_2 )$  be two groupoids over $X$.
A groupoid morphism {\it over } $X$ or $X$-{\it morphism}
\[ \phi \colon G_1 \rightarrow G_2  \]
is a functor such that $Ob(\phi ) = 1_X$. }
\end{defn}
\begin{defn} {\rm
The category of groupoids, denoted by $Grd$, has as its objects
all groupoids
$G = (G, X, \alpha, \beta, m, i)$ and its morphisms are functors between them.

For each space $X$, the category of groupoids over $X$, denoted by $Grd(X)$,
has as its objects groupoids with object set $X$ and $X$-morphisms as its morphisms.}
\end{defn}
\begin{defn}{\rm
A {\it topological groupoid} is a groupoid $G = ( G, X, \alpha,
\beta, m, i )$ where $G$ and $X$ are sets equipped with a
topology, such that the four structure maps and the inverse map
$G\rightarrow G$, $g\mapsto g^{-1}$, assigning to each arrow $g$
its inverse are continuous with respect to these two topologies. A
morphism of topological groupoids is a continuous functor between
them.}
\end{defn}
\begin{defn}{\em A topological groupoid $G$ on $X$ is called {\em \'etale}
if the source $\alpha :G\to X$ is a local homeomorphism (This implies
that all other structure maps are local homeomorphisms, also).}
\end{defn}
For more information and examples of \'etale groupoids,
we can give Moerdijk' papers such as, see \cite{Mo1,Mo2,Mo3}
\begin{defn}{\rm
A topological groupoid $G = (G, X, \alpha, \beta, m, i)$  is said to be
{\it open} if the source and target maps $\alpha, \beta$ are both open
maps. }
\end{defn}

We have a category $TGrd$ whose objects are topological groupoids and
morphisms are continuous functors.

For each topological space $X$, the category of topological
groupoids over $X$, denoted by
$TGrd(X)$ or $TGrd/X$, has its objects topological groupoid with object $X$
and continuous $X$-morphism as its morphisms.

Notice that the null topological groupoid $X$ is an initial object
and also the topological groupoid $X\times X$ which is called coarse
topological groupoid is a final object in the category $TGrd(X)$.

We can give some examples of topological groupoids mentioned earlier.
In Example \ref{143}, we consider the set $X$ to be a topological space,
then $X$ can be regarded as a topological(\'etale) groupoid.
Also in Example \ref{145}, the groupoid $X\times X$ is a topological groupoid on a
topological space $X$. In addition, an equivalence relation $R$ is a topological
groupoid over the topological space $X$.
The fundamental groupoid $\pi_1 (X)$ given
in Example \ref{147} is a
topological groupoid in a natural way given certain
local conditions on a topological space $X$ \cite{Br-Da}.
Also each topological group provides an example of a topological
groupoid.
 If $ * \colon G\times X\rightarrow X$ is a continuous action of
topological group $G$ on a topological space written $*(g, x) = gx$,
then the product $G\times X$  has a topological groupoid structure
on $X$
\cite{Br3} with the multiplication $(h, gx)(g, x)= (hg, x)$. Its source and target maps are the second projection and
the action itself, respectively.

%\end{document}

\chapter{Local and Global Subgroupoids}

In this chapter, we give the definition of a local subgroupoid of
a groupoid $G$ on  a topological space $X$. We show that many
ideas generalise  smoothly from the local equivalence relations to
the local subgroupoids (recently, the concept has been explored
widely in \cite{Br-Ic,BIM}). Now here we shall give some of the
relations between the local subgroupoids and the local equivalence
relations.

Firstly, we give the definition of a local subgroupoid and
 some main examples.
\subsection{Local subgroupoids}
Let $X$  be a topological  space  and  let $G$ be  a  groupoid  with
$Ob(G)$ = $X$. Let $U$ be an open subset of $ X$. Let  $L_G(U)$
be the set of all wide
subgroupoids of $ G|U$, where $G|U$ is the full
subgroupoid of $ G$ on $U$, i.e.,
$G|U = \alpha^{-1}(U)\cap \beta^{-1}(U)$.
Let $V, U\subseteq X$ be open sets such that $V\subseteq U$.
If $H$ is a wide subgroupoid of $G|U$, then $H|V$ is a wide subgroupoid
of $G|V$,
So there is a restriction map
\[
\begin{array}{ccc}
        L_{UV}\colon L_G(U) &\longrightarrow &L_G(V) \\
                                         H &\longmapsto        & H\mid V.
\end{array} \]
and these define a presheaf
\[L_G \colon {\cal O}(X)^{op}\longrightarrow Sets. \]
However $ L_G$ is not in general a sheaf, as explained in
Chapter 1, see also \cite{BIM}.
In the usual way this presheaf $L_G$ defines a sheaf denoted by ${\cal L}_G$.
We can define ${\cal L}_G$ as follows:
\[     {\cal L}_G = \bigcup_{x\in X}{{\cal L}^G}_x  = \bigcup_{x\in X}
 \{ (U_i, H_i)_x \colon x\in U_i\subseteq  X \mbox{ open }  \ H_i \in L_G(U_i) \}  \]
                           and
\[  p\colon {\cal L}_G\rightarrow X,   \  \  \ p({\cal L}_x) = x        \]
is the canonical projection, i.e., it is a local homeomorphism.
Let $U\subseteq X$ be open and $s\colon U\rightarrow {\cal L}_G $ be a section
of ${\cal L}_G$ over $U$. The set of such sections is denoted by
$\Gamma (U, {\cal L}_G)$ which defines a presheaf:
\[   \Gamma \colon {\cal O}(X)^{op}\longrightarrow Sets .  \]
The set of global sections of ${\cal L}_G$ is denoted by $\Gamma (X,{\cal L}_G)$
in Chapter 1.
Also every element $H\in
L_G(U)$ is associated with a section $s\in \Gamma(U,{\cal L}_G)$. If
$x\in X$ and $\sigma \in{\cal L}_x $, then there are an open
neighbourhood  $x\in U\subseteq X$ and an $s \in \Gamma (U,{\cal L}_G)$ such that
\[      \sigma = (U, H)_x = germ_xH = s(x)=s_x,      \]
for more detail, see Chapter 1.
\begin{defn}{\rm
A {\it local subgroupoid} of a groupoid $G$ on the topological space $X$ is a global  section  of  the
sheaf ${\cal L}_G$ associated to the presheaf $L_G$.}
\end{defn}

In other word, a local subgroupoid of $G$ is a continuous functions from $X$ to
${\cal L}_G$ such  that $p\circ s =id_X$.

A local subgroupoid $s$ can be defined by atlas as follows.
Given an open cover $\{U_i:i\in I\}$ of $X$,
subgroupoids $H_i$, $H_j$ on $U_i$, $U_j$, respectively,
such that each point $x\in U_i\cap U_j$
has a neighbourhood $W$ on which $H_i$ and $H_j$ agree .
This
atlas is  denoted by ${\cal U}_s=\{(U_i,H_i):i\in I\}$.
We define an equivalence relations on the atlas, namely,
$H_i$ and $H_j$ are equivalent
at $x$ if there exists $W$ such that $x\in W\subseteq U_i\cap U_j$  and $H_i$ and $H_i$  agree on $W$.
The equivalence classes are called the {\em germs} of
subgroupoid of $G$  at $x$. These germs form the stalk ${{\cal L}^G}_x$.

Now  we give examples of local subgroupoids.
\subsection{Examples}

\begin{example}{\rm
Let $X$ be a topological space.
Every local equivalence relation on $X$ is a local subgroupoid
of  the groupoid $X\times X$.
 For open set $U$ in $X$, let $E(U)$ be the set all equivalence
relations on $U$.
This $E : {\cal O}(X)^{op} \rightarrow Sets $ is a presheaf which
defines a sheaf ${\cal E}$ on $X$.
So we obtain a local equivalence relation $r$ as a global section of ${\cal E}$.
It is easy to show that the set $E(U)$ is the family of all wide
subgroupoid of $U\times U$. Hence a local equivalence relation
$r$ on $X$ is a local subgroupoid of
$X\times X$.}
\end{example}
\begin{example}\label{223}
{\em Any topological space $X$ can be considered as a groupoid on
itself with $\alpha = \beta = id_X$ every element the identity,
see Example \ref{143}, called null groupoid.  As is well-known, it
is an initial object in the category of groupoid on $X$, $Grd(X)$.

Let us build up the local subgroupoid of $X$ on $X$. Firstly,
we shall construct the corresponding sheaf ${\cal L}_X$. Let $U\subseteq X$
be open, the set
\[ L_X(U) = \{ U : U \mbox{ \ \ is a wide subgroupoids of} \ \   X|U = U \} \]
is just restriction set of $U$; $L_X(U) = \{ U \}$.
For open sets $V, U \subseteq X$ with $V\subseteq U$, the restriction map is
\[  L^X_{UV} \colon  L_X(U)\rightarrow   L_X(V), \ \ \  U\mapsto V,   \]
i.e., $L^X_{UV}(U) = U|V = V$. Hence
$L_X\colon {\cal O}(X)^{op}\rightarrow Sets $ is a presheaf.
Let ${\cal L}_X$
be the sheaf associated with presheaf $L_X$. Its stalks are such as
\[   {{\cal L} ^X}_x = \lim_{\stackrel {\longleftarrow }{x\in X}}L(U) =
\lim_{\stackrel {\longleftarrow}{x\in U}} \{U\}= U              \]
and the sheaf
\[ {\cal L}_X = \bigcup_{x\in X} {{\cal L}^X}_x =
\bigcup_{x\in X} \{(U, U)_x : x\in U, U\in L_X(U)    \}   \]
So a local subgroupoid of $X$ is defined  by an atlas
${\cal U}_s=\{(U_x, U_x): x\in X\}$.}
\end{example}
The following examples give us the relations between local subgroupoids
and local equivalence relations.

\begin{example}{\em
Let $X^* = \{ X_i \mid i\in I \}$ be a
partition of a topological space $X$ and let $R$ be  an associated
equivalence relation of partition
on $X$ . Let $K$ be a group, then $G = R\times K$ is
a groupoid on $X$, see Example \ref{group} .

Let $U_i, U_j\subseteq X$ be open sets with ${U_j}\subseteq {U_i}$.
Let us consider the
following structure on $X$. We define the sets
\[  E(U_i) = \{ R_i \mid R_i \mbox { is an equivalence relation on } U_i \} \]
and similarly
\[  L(U_i) = \{ H \mid H \mbox { is a wide subgroupoid of}  \  G\mid_{U_i} =
R\times K\mid_{U_i} = R\mid_{U_i}\times K \}.    \]
Clearly, they define functors from ${\cal O}^{op}$ to $Sets$, i.e., $E$ and
$L$ are presheaves on $X$. So  we define a natural transformation
$\gamma \colon E\rightarrow L$,
by $\gamma_U\colon E(U_i)\rightarrow L(U_i)$, $R_i\mapsto R_i\times K$,
for each open set $U_i\subseteq X$. In other word, the following diagram is
commutative;
$$
\begin{diagram}
               E(U_i)& \rArr^{\gamma_{U_i}}     &L(U)\\
        \dArr^{{\cal E}_{U_iU_j}} &                 &\dArr_{L_{U_iU_j}} \\
             E(U_j)  & \rArr_{\gamma_{U_j}}     &L(U_j)
\end{diagram}
\begin{diagram}
               R_i      &     \rMapsto    & R_i\times K \\
               \dMapsto &                 & \dMapsto    \\
               R_i\mid_{U_j}&\rMapsto     & R_i\mid_{U_j}\times K = R_i \times K\mid_{U_j}
\end{diagram}
$$
That means that $\gamma \colon E\rightarrow L$ is a presheaf morphism. It give
rise to a sheaf morphism between associated  sheaves \cite{Br}. Let ${\cal E}_R$ and
${\cal L}_G$ be the sheaves associated to the presheaves  $E_R$ and $L_G$, respectively.
Let $\gamma^* \colon {\cal E}_R \rightarrow {\cal L}_G$ be
the corresponding sheaf morphism, i.e., ${\gamma^*}$  is continuous and the following diagram
$$
\begin{diagram}
      {\cal E}_R    &                 &\rArr^{\gamma^*}&     & {\cal L}_G\\
                    & \SE_{p_{\cal E}} &                &  \SW_{p_{\cal L}} & \\
                    &                 &         X        &     &
\end{diagram}
$$
is commutative.

If $r$ is a local equivalence relation of $R$ given by an atlas
${\cal U}_r = \{ (U_i, R_i) \mid i\in I\}$, then $\gamma^* r$ is to be a
local subgroupoid of $G$ on $X$. It can be defined as  an atlas
${{\cal U}_r}^* = \{ (U_i, R_i\times K) \mid i\in I \}$. In other word, a
local equivalence relation of  $R$ on a topological
space $X$ defines a local subgroupoid of groupoid   $R\times K$ on $X$,
while $K$ is a group.}
\end{example}

\begin{example} \label{225}{\em
Let $G$ be a groupoid on a topological space $X$.
Then $X\times X$ is also a groupoid on $X$, see Example \ref{145}.

For $U\subseteq X$ open set, we obtain the following presheaves.

\[       L_G = \{L_G(U), L_{UV},  X\}  \ \ \ \mbox{and} \ \
E_{X\times X} = \{ E(U), E_{UV}, X \}
\]
respectively. The set $L_G(U)$ is all wide subgroupoids of $G|U$,
$E(U)$ is all equivalence relations on $U$.
These presheaves define the sheaves on $X$.
We have defined a local
equivalence relation $r$ on $X$ to be a global section of the sheaf ${\cal E}_{X\times X}$
associated with the presheaf $E_{X\times X}$.

Likewise a local subgroupoid  $s$ of $G$ is a global section of the sheaf ${\cal L}_G$
associated with $L_G$. Also we can define a morphism of groupoids on $X$ as follows;
\[ \Upsilon = [\alpha, \beta] \colon G \longrightarrow X\times X     \ \ \
g\mapsto  (\alpha(g), \beta(g))   \]
where $\alpha $ and $\beta $ are  the source and target maps, respectively \cite{Mac}.
Let $H\in L(U)$ and $R\in E(U)$  such that $ \Upsilon [\alpha,\beta](H)=R$.
So we obtain a morphism of groupoids $H\rightarrow R$ on $U$ by $\Upsilon =[\alpha, \beta]$.
Hence we define a map $\Upsilon \colon L(U)\rightarrow E(U)$, $H\mapsto R$,
 which defines a morphism of presheaves which is a natural transformation of functor.
In fact, the diagram
$$
\begin{diagram}
               L_G(U)& \rArr^{\Upsilon_U}     &E(U)\\
        \dArr^{{\cal E}_{UV}} &                 &\dArr_{E_{UV}} \\
             L_G(V)  & \rArr_{\Upsilon_V}     &E(V)
\end{diagram}
$$
is commutative. Then the natural transformation
\[
\Upsilon \colon L_G\longrightarrow E_{X\times X}
\]

defines a sheaf morphism
\[
\Upsilon^*\colon{\cal L}_G\longrightarrow {\cal E}_{X\times X}.
\]
i.e., $p_{X\times X}\circ \Upsilon^* = p_G$, where $p_{X\times X}$ and
$p_G$ are local homeomorphism on sheaves ${\cal E}_{X\times X}$,
${\cal L}_G$, respectively. So we obtain $\Upsilon^* (s) =r$}
\end{example}

In this example, the idea leads us to a  generalization.
Namely, any groupoid morphism gives rise to a sheaf
morphism between the corresponding sheaves. In other word,
there is a functor from the
category of groupoids $Grd(X)$ to  the category of sheaves $Sh(X)$.

Suppose given two groupoids $G, H$ and a groupoid morphism
$\Phi \colon G\rightarrow H$ in $Grd(X)$. As usual, we can obtain
the following presheaves of sets on $X$ for the groupoids $G$ and $H$:
\[    L_G = \{ L_G(U), {L^G}_{UV}, X \} , \ \ \ L_H = \{ L_H(U), {L^H}_{UV}, X \}        \]
for open sets $U, V\subseteq X$ with $V\subseteq U$.
That is, $L_G$ and $L_H$ are functors from
${\cal O}^{op}$ to $Sets$. Hence the morphism of groupoid
$\Phi \colon G\rightarrow H$ gives rise to a morphism of presheaves and
a natural transformation by the following diagram

$$
\begin{diagram}
               L_G(U)&   \rArr^{{\Phi^*}_U}     &L_H(U)\\
   \dArr^{{L^G}_{UV}} &                          &\dArr_{{L^H}_{UV}} \\
             L_G(V)  &    \rArr_{{\Phi^*}_V}     &L_H(V)
\end{diagram}
$$ The natural transformation $\Phi^* \colon L_G\rightarrow L_H$
is explicitly  defined as follows. For each open set $U\subseteq
X$, suppose that $A\in L_G(U)$, i.e., if $A$ is a wide subgroupoid
of $G|U$, then its image $\Phi_U(A)$ is a wide subgroupoid of
$H|U$, since $Ob(\Phi^*)$ is the identity on $U$. Thus
$\Phi_U(A)\in L_H(U)$. Finally we have $\Phi^*(A|U) =
\Phi_U(A)|V$.

In sheaf theory, as is well-known,
a morphism of presheaves defines a morphism of corresponding sheaves
\cite{Br,Sw,Ma-Mo}.

Let ${\cal L}_G$ , ${\cal L}_H$ be the sheaves associated with presheaves
$L_G$ and $L_H$, respectively , and also let $\Phi^* \colon L_G\rightarrow L_H$
be the presheaf morphism as above. Then for each $x\in X$,
$\Phi^*$ induces a morphism
\[  {\Phi^*}_x\colon {{\cal L}^G}_x = \lim_{\stackrel {\longrightarrow}{x\in U}}
L_G(U) \rightarrow    \lim_{\stackrel {\longrightarrow}{x\in U}}
L_H(U) = {{\cal L}^H}_x        \]
\[   (U, A)_x\rightarrow ((U, \Phi (A))_x                  \]
and therefore a map $\Phi^* \colon {\cal L}_G\rightarrow {\cal L}_H$
which is a sheaf morphism.

We summarise all knowledge above as follows:
Given two groupoids $G$,  $H$  and a groupoid morphism $\phi \colon G\rightarrow
H$ in the category of groupoids on $X$, $Grd(X)$. Then morphism give rise to  a functor
${\cal L}$ from $Grd(X)$ to the category of sheaves $Sh(X)$.

Also we obtain a subcategory of $Sh(X)$ whose objects are sheaves
defined by groupoids as  above
and whose  morphisms are sheaf morphism. It is denoted by $Sh_{Grd}(X)$ and is
a full subcategory of the category $Sh(X)$.

\begin{prop}
Let $Sh_{Grd}(X)$ be the subcategory of the category $Sh(X)$ as above.
Then its initial and final objects are ${\cal L}_X$ and ${\cal L}_{X\times X}$,
respectively.
\end{prop}
\begin{pf}
First we shall show that the sheaf ${\cal L}_X$ is
an initial object of the category  $Sh_{Grd}(X)$.
Clearly, ${\cal L}_X$ is an object of the category  $Sh_{Grd}(X)$. Now we
have to show that for each object ${\cal L}_G$ in   $Sh_{Grd}(X)$, there
is a unique sheaf morphism between sheaves  ${\cal L}_X$ and ${\cal L}_G$. That is,
there is a unique map $\phi \colon {\cal L}_X\rightarrow {\cal L}_G$
such that the following diagram
$$
\begin{diagram}
              {\cal L}_X  &           &\rArr^{\phi^*} &   &  {\cal L}_G   \\
                          & \SE_{p_X} &                 &\SW_{P_G} &        \\
                          &           &      X           &           &
\end{diagram}
$$
commutes.

Since $X$ is an initial object of the category $Grd(X)$, there is
a unique groupoid morphism  $I\colon X\rightarrow G$, for each
groupoid $G$ in $Grd(X)$. This morphism give rise to a sheaf
morphism $I^*\colon {\cal L}_X\rightarrow {\cal L}_G$ with  $ p_G
o I^* = p_X$. Therefore $I^*$ must be unique, since $I$ is unique.
So $I^* = \phi^*$. Here  the sheaf morphism  ${I^*}$ is defined as
follows, on each stalks ${{\cal L}^X}_x$ ${I^*}_x\colon {{\cal
L}_x}^X\rightarrow {{\cal L}_x}^G$, $(U, U)_x\mapsto (U, I(U))_x$.

Likewise the sheaf ${\cal L}_{X\times X}$ can be showed to be a final object of
the category $Sh_{Grd}(X)$.
\end{pf}

{\bf Category of Local Subgroupoids : LSG(X)}

Let $Sh_{Grd}(X)$ be the category of sheaves defined by groupoids  as above.
We can define a category of local subgroupoids and denote by $LSG(X)$
in which an arrow $\phi^* \colon s_1\rightarrow s_2$ is a sheaf morphism
$\phi^*\colon {\cal L}_G\rightarrow {\cal L}_H$ with $s_1 o \phi^* =
s_2$ in the category $Sh_{Grd}(X)$, i.e., the following diagram
$$
\begin{diagram}
             {\cal L}_G &           & \rArr^{\phi^*}&  & {\cal L}_H \\
                        &\NW_{s_1}  &  &              \NE_{s_2}& \\
                        &           &       X        &  &
\end{diagram}
$$
commutes, for each groupoid $G$ and $H$ in the category $Grd(X)$.
Its class of
objects is the set of local subgroupoids of $Grd(X)$.
%---------------------------------------------------------------------------
\section{Coherent Local Subgroupoids}

In this section, we generalise  the idea of coherence from local
equivalence
relations, as  given by Rosenthal \cite{Ro2},
to local subgroupoids.

We  define a partial order
structure on the sheaf of germs  of  the local subgroupoids as follows.

Let $s, t$ be  local subgroupoids of $G$, and let $s_x,\ t_x$ be their germs at $x\in X$.
Then there are open sets $U_x$ and $V_x$  containing $x$ and
$H_x\in L_G(U_x), \  K_x\in L_G(V_x)$ such that $H_x$ defines $s_x$
and $K_x$ defines $t_x$, i.e.
\[       s_x =(U_x, H_x)_x \  \mbox{ and}   \    t_x = (V_x, K_x)_x.
\]
We say  $s_x\leq t_x$ if there is an open neighbourhood $W$ of $x$ with
$W\subseteq U_x\cap V_x$ and  $H_x\mid W \subseteq K_x\mid W$, i.e. $H_x$
is a wide subgroupoid of $K_x$ on $W$. So we obtain a natural order structure on the
sheaf of germs of local subgroupoids by the following definition.
\begin{defn}{\rm
Let $s$ and $t$ be local subgroupoids of $G$ and let  $s = (s_x)_{x\in X}$
and $t = (t_x)_{x\in X}$. We write $s\leq t$ iff  $s_x\leq t_x$ for all $x\in X$. }
\end{defn}
\begin{defn}{\rm
Let $L_G(X)$ be the set of wide  subgroupoids of $G$ on $X$ and let
$H\in L_G(X)$ and $x\in X$. Then  $ loc(H)$
is the local subgroupoid defined by
\[             loc(H)(x) = (X, H)_x.    \]}
\end{defn}
 For  the open set $U$ of $X$, clearly
$loc(H)(U) = (U, H|U)$
 \cite{BIM}.

Let $s$ be a local subgroupoid of  $G$ on $X$. Globalisation of  $s$ is
defined
by
\[    glob(s) = \cap \{ H : s\leq loc(H) \}          \]
where $H\in L_G(X)$. Thus  $s\leq loc(H)$ means that for $x\in X$,
 $s_x\leq (loc(H))_x$, i.e., if $s_x = (U_x, H_x)_{x\in X}$,
\[
       H_x|{U_x} = H_x \subseteq (loc(H))_x| {U_x} = H|{U_x} .
\]
We  think of glob(s) which obtains approximate $s$ by a global subgroupoid.
To get a best possible approximate we consider all wide subgroupoids of $G$
on $X$  which  locally contain $s$, and intersect them. That is, $glob(s)$
is the smallest wide subgroupoid of $G$ on $X$ which locally contains  $s$.

Moreover we can provide an alternative useful description of glob(s).
Let us suppose that $s$ is given by an atlas
${\cal U}_s=\{(U_x, H_x): x\in X\}$, i.e.,
\[           s = (s_x)_{x\in X} = (U_x, H_x)_{x\in X}.   \]
Let $ {\cal V} = \{V_x : x\in X\}$ be an open cover of $X$ such that
$x\in V_x\subseteq U_x$ for $x\in X$.  Let $H_{\cal V}$ be the subgroupoid
of $G$ generated by   $\{ H_x|{V_x} : x\in X \}$. Clearly,
$H_{\cal V}\subseteq G$,
 because  $H_x|{U_x} = H_x \supseteq H_x|{V_x}$ and then
\[
      \bigcup_{x\in X}H_x|{V_x} =
      \bigcup_{x\in X}H_x = H_{\cal V}\supseteq G.
\]
\begin{example}{\rm
Let $X$ be a topological space and let
${\cal U} =\{ U_x : x\in X \}$ be an open cover of $X$.
We consider the groupoid $G = X\times X$ on $X$. Then
\[   G|{U_x} = U_x\times U_x ={\pi_1}^{-1}(U_x)\cap {\pi_2}^{-1}(U_x)      \]
where  $\pi_1$ and $\pi_2$ are the natural projections. But the subgroupoid
$H_{\cal U}$ is generated by $\{G|{U_x}\}$, i.e.,
\[     H_{\cal U} = \bigcup_{x\in X}G\mid_{U_x} = \bigcup_{x\in X}(U_x\times U_x)  \]
and  $H_{\cal U}\subset G$.}

\end{example}

\begin{prop}
Let $s$ be a local subgroupoid of $G$ given by an atlas
${\cal U}_s=\{(U_x, H_x): x\in X\}$. Let ${\cal U}=\{U_x:x\in X\}$ and
 let $H_V$ be the subgroupoid of $G$ generated by $\{H_x|{V_x}:x\in X\}$.
 Then
\[  glob(s) = \bigcap \{ H_{\cal V }: {\cal V}\leq {\cal U} \}  \]
\end{prop}
\begin{pf}
Let $Q$ be a subgroupoid of $G$ on $X$ such that $s\leq loc(Q)$. Take
an open cover $\{ W_x : x\in X \}$ with $x\in W_x\subseteq U_x$ and so
$H_x|{W_x}\subseteq Q|{W_x}$. Then $H_W\subseteq Q$ and hence
$\bigcap \{H_V :{\cal V}\leq {\cal U} \}\subseteq glob(s)$.

Conversely, if ${\cal V}\subseteq {\cal U}$, then
$H_x|{V_x}\subseteq H_{\cal V}|{V_x}$,
since $H_{\cal V}$ is locally generated by $\{ H_x|{V_x} \}$. Thus
$s\leq loc(H_{\cal V})$. So
\[               glob(s) = \bigcap \{H_{\cal V} : {\cal V}\leq {\cal U} \}            \]
\end{pf}

We always have  $glob(loc(H))\subseteq H$, for each $H\in L_G(X)$.
In fact, we have $loc(H)\subseteq H$ , and  it follows that
$glob(loc(H))\subseteq glob(H)=H$
\cite{Ro1, BIM}.

\begin{defn}\label{235}
Let $s$ be a local subgroupoid of $G$ on $X$.

i)  s is called {\it coherent} if $s\leq loc(glob(s))$.

ii) s is called {\it globally coherent} if $s = loc(glob(s))$

iii) s is called {\it totally coherent} if for every open set $U, s\mid_U$
is coherent.
\end{defn}
\begin{defn}\label{236}
Let $H$ be a subgroupoid of $G$ on $X$, i.e. $H\in L_G(X)$.

i)  $H$ is called {\it locally coherent} if $loc(H)$ is coherent.

ii) $H$ is called {\it coherent} if $H = glob(loc(H))$.
\end{defn}

Coherence of $s$ says that in passing between local and global information
nothing is lost due to collapsing. Note that for a groupoid $H$  to be coherent
we must have that for every open cover ${\cal V} = \{V_x : x\in X \}$, $H =H_{\cal V}$, where $H_{\cal V}$ is
the subgroupoid of $H$ generated by $\{ H\mid_{V_x}:x\in X \}$. We can find many examples for local and
global case in  Rosenthal's papers \cite{Ro1, Ro2}, if we take a local
equivalence relation as a local subgroupoid.
(In the local groupoid case
have been considered  widely in  paper  `Local subgroupoids II: examples
and properties \cite{BIM}.)

\begin{example}{\rm
Any topological space $X$ can be considered as
a groupoid on itself \ref{143}.
Let ${\cal L}_X$ be a sheaf corresponding to the groupoid $X$,
 see Example  \ref{223}.
Let $s$ be a local subgroupoid of $X$.
 It is  easily seen  that  $loc(X) = s$, that is,
 $loc(X)(U)=(U, X|U) = (U, U)$.
Moreover $glob(s) = X$, since $glob (s) = glob (loc (X)) = X$. So
$s$ is globally coherent and $X$ is coherent.}
\end{example}

A bundle of groups is a good example of a groupoid.
A  bundle  of  groups  can  also  described as  a  bundle   $p\colon
G\rightarrow X$ of
spaces in which each fiber $p^{-1}(x)$ is a group in such a way that the resulting
operations of addition $+\colon G\times_p G\rightarrow G$ and inverse
$-\colon G\rightarrow G $ are continuous maps \cite{Ma-Mo}.
\begin{example}{\rm
Clearly a bundle of groups is a groupoid  whose source and target
maps are equal, i.e. $\alpha = \beta = p$.
Let $U$ be an open set in $X$. The set
\[  p^{-1}(U) = G|U = \bigcup_{x\in U}p^{-1}_x = \bigcup_{x\in U} G_x \]
is a groupoid on $U$. For each open set $U$ in $X$, we obtain $L_G(U)$
of the set of all subbundles of groups of $G|U$ on $U$.
For $V\subseteq U$ open sets in $X$,
\[
\begin{array}{ccc}
        L_{UV}\colon L_G(U) &\longrightarrow &L_G(V) \\
                                        (G|U) &\longmapsto & (G|U)|V
\end{array} \]
is a restriction morphism which  defines a presheaf
\[   L_G : {\cal O}(U)^{op}\rightarrow  Sets.              \]
 Then we get a sheaf ${\cal L}_G$ from the
presheaf $L$.
Let $s$ be a local subgroupoid of the bundle of groups $G$ on $X$.
Since $L_G(U) =\{ G|U\}$, then $s$ is a globally coherent local
subgroupoid of  $G$ and $G$ is a coherent groupoid on $X$. In fact,
now, let ${\cal V} = \{ V_x  : x\in X\}$ be an open cover
of $X$ such that for each $x\in X$, $x\in V_x\subseteq U_x$, where
$ {\cal U} = \{ U_x : x\in X \}$ is also open cover  of $X$.
Let $H_{\cal V}$ be the subgroupoid of $G$ generated
by $\{G|V_x : x\in X\}$.  $p^{-1} (U_x) = G|{U_x}\supseteq H_x|{V_x}$.
But $H_{\cal V} = G$ and  $glob(s) = \cap \{ H_V : V\leq U \} = H_{\cal V} = G$.
So $loc(glob(s)) = loc (G) = s$. Hence $s$ is a globally coherent. Since
$loc(G) = s$ and $glob(loc(G)) = glob(s) = G$.
}
\end{example}
(More examples  of local subgroupoids are given  in \cite{BIM}.)

We obtain functors loc and glob as follows:  Let {\bf CL} be the category of coherent local
subgroupoid and {\bf CG} be the category of locally coherent global subgroupoids
on $X$.
\begin{prop}
Let  $glob\colon {\bf CL}\longrightarrow {\bf CG}$ and
$loc\colon {\bf CG}\longrightarrow {\bf CL}$ be functors. These functors
form a pair of adjoint functors with $loc$  left adjoint to $glob$, i.e.
$ loc\dashv glob$.
\end{prop}
\begin{pf}
Our categories are  ${\bf CL} = \{ s : s\leq loc(glob(s))\}$ and
${\bf CG} = \{ H : loc(H) \mbox{is coherent} \}$, and natural bijection
\[   \theta\colon {\bf CL}(loc(H), s)\longrightarrow {\bf CG}(H, glob(s))   \]
where $H\in {\bf CG}$ and $s\in {\bf CL}$. If we take $H = glob(s)$,
it gives a unique map
\[   \eta \colon loc(glob(s) \longrightarrow s       \]
such that $\theta(\eta) = I$. This map $\eta$  is a unit of adjunction such that
$s\leq loc(glob(s))$ and similarly for $loc(H)$, gives a unique map
\[   \xi \colon H\longrightarrow glob(loc(H))      \]
The map is dual to the unit of an adjunction is the counit such that
$glob(loc(H))\subseteq H$
\end{pf}

Under this adjunction, we have an equivalence between globally coherent
$s$ and coherent $H$.

\begin{cor}\label{239}
Let $H\in L_G(X)$ and $ s = loc(H)$. Then the transitivity components
of $H_{\cal V}$ are relatively closed and open in the transitivity
components of $H$.
\end{cor}
\begin{pf}
Let $M_{x,v}$ and $M_x$ denote the transitivity components of $x$ in $H_{\cal V}$
and $H$ respectively. Clearly,  $M_{x,v}\subseteq M_x$. Because,
$s = loc(H)$, $glob(s) = glob (loc(H))\subseteq H$, so $H_{\cal V}\subset H$.
Let $y \in M_{x,v}$. Then there is $x_1,x_2,...,x_n,x.$, \ $V_1,V_2,...,V_{n+1} $
such that  $h_1\in H\mid_{V_1}(y,x_1),  h_2\in H\mid_{V_2}(x_1, x_2),...,
 h_n\in H\mid_{V_{n+1}}(x_n, x)$, i.e.  $g = h_n...h_2.h_1$.
Take $V_1\cap M_x$ and let $z\in V_1\cap M_x$. Hence $h\in H(z,x)$ and since
$k\in H(y,x)$, $h^{-1} k \in H_V(y, z)$. Thus, $h^{-1} k\in H_{\cal V}$ and $z\in M_{x,v}$.
 Now let us show that $M_{x,v}$ is closed in $M$. Let $z\in\overline{M_{x,v}}$
be the closure relative to $M_x$. For every open neighbourhood $U$ of $x$, we have
$V_z$,  take an element $y\in V_z\cap M_{x,v}$. Then, there is a $g\in H_v(y,x)$ and since
$h\in H(z, x)$, we have $h^{-1} g\in H(y,z)$. Since $y,z\in V_z$,
 $h,g\in H_V(y,z)$ and $z\in M_{x,v} $. Thus $M_{x,v} = \overline{M_{x,v}}$
\end{pf}
\begin{thm}\label{2310}
Let $H\in L_G(X)$ and $ s = loc(H)$. If
there is an open neighbourhood $W_x$ of $x\in X$ such that $H|{W_x}$ has
connected transitivity
components, then $s$ is coherent.
\end{thm}
\begin{pf}
Suppose that $s$ is not coherent. Then, for some $a\in X$, we have
$s_a\not\leq loc(glob(s))_a$, i.e. given any open neighbourhood $W$ of $a$,
there is a cover $\{ V_x :x\in X\}$  and $y, z\in W$ such that there
exists an $h\in H(y,z)$, $h\not \in H_{\cal V}(y,z)$. In particular, this is true for $W_a$.
By Corollary \ref{239}, the transitivity component of $y$ in $H\mid_{W_a}$ is clopen
in the transitivity component of $y$ in $H|{W_a}$, which is connected.
This is a contradiction as it forces  $h\in H(y,z)$.
\end{pf}
\begin{cor}
\label{2311}
Let s be a local subgroupoid of G defined by $H_x\in L_G(U_x)$, $x\in X$. Suppose that for
every $a\in X$ , every open $V\subseteq U_a$ with $a\in V$, there is an open
neighbourhood W of $a$ with $W\subseteq V$ and with $ H_a\mid_V$ having connected transitivity
component. Then $s$ is totally coherent, i.e. $s\mid_U$ is coherent for
every open  U in X.
\end{cor}
\begin{pf}
If $a\in U$, where $U$ is open, consider $H_a\mid_{U\cap U_a}$ and apply the argument from
the Theorem \ref{2310}.
\end{pf}
\begin{thm}\label{2312}
Let $H\in L_G(X)$ and suppose $H$ has connected transitivity components.
Then $H = glob(loc(H))$. Conversely, if $H = glob(loc(H))$ and $H$ has closed transitivity
components, then it has connected transitivity components.
\end{thm}
\begin{pf}
Given an open cover ${\cal V}=\{ V_x :x\in X\}$ of $X$. The subgroupoid $H_{\cal V}$
generated by  $\{H|{V_x} \}$ is contained in $H$, $H_{\cal V}\subseteq H$
and,  by Corollary \ref{239}, since transitivity components  of  $H_{\cal V}$
are relatively in those of $H$,
which are connected, since  $M_x$  is connected, it must be $M_{x,v} = M$
so we have must $H_{\cal V} = H$,  $H = glob(loc(H)) = \cap H_{\cal V}$.

If $H = glob(loc(H))$, for every cover ${\cal V}$, $H_{\cal V} = H$ .
Let $a\in X$ be such that $M_a$, the transitivity component of $a$ in $H$, is
not connected. Let $U$ and $V$  be open sets separating $M_a$.
Let ${\cal U} = \{ U, V, X-\{x\} \}$. Choose $x,y \in M_a$ such that
$x\in U\cap M_a$,  $y\in V\cap M_a$. Then there exists $g\in H(x,y)$ but
$g\notin H_{\cal V}(x,y)$ since $(U\cap M_a)\cup (V\cap M_a) = M_a$ and
they are disjoint, since $H_{\cal V}\not\subseteq H$ we have that  $glob(loc(H))\subseteq H$.
This is a contradiction.
\end{pf}
\begin{prop}
i) Suppose $s$ is globally and totally coherent on $X$. If $U$ is open
in $X$, then $s|U$ is globally coherent.

ii) If there is an open cover $\{ V_x : x\in X \}$ of $X$ such that
$s|{V_x}$ is globally and totally coherent for all $x\in X$, then $s$ is totally coherent.
\end{prop}
\begin{pf}
i) We  have  $s = loc(glob(s))$. By definition,
$glob(s|U)\subseteq glob(s)|U$,
hence $loc(glob(s|U))\leq loc(glob(s)|U) = loc(glob(s))|U = s|U$.
Sine $s\mid_U$ is coherent, by totally coherence of $s$, we have
$s|U\leq loc(glob(s|U))$. So $s|U = loc(glob(s|U))$, i.e.
$s|U$ is globally coherence.

ii) By (i), if $U$ is open in $X$ and $s|{V_x}$ is globally
coherence for all
$x\in X$, then $s|{U\cap V_x}$ is globally coherent. Thus  $s|{U\cap V_x}$
$= loc(glob(s))|{U\cap V_x} \leq loc(glob(s|U))|{V_x}$,
since this holds for all $x\in X$,
we have $ s|U\subseteq loc(glob(s|U))$, i.e. $s$ is totally coherent.
\end{pf}

\subsection{ Topological foliations}
One of Ehresmann's approaches to the foundations of foliation
theory goes via the consideration of a topological space equipped
with a further `fine' topology. Such fine topologies appear also
in the context of local equivalence relations and have been
considered in \cite{AGV} and in \cite{Ro1}.We shall need the
following elaboration of this idea for the local subgroupoids.

Let $s$ be a local subgroupoid of $G$ on a topological space $X$
which is given by an atlas $\{(U_x, H_x): x\in X\}$. We can define
a new topology on $X$ denoted by $X^s$. The underlying set of
$X^s$ is $X$. Let $M_{x,a}$ denote the transitivity components of
$x$  in $U_a$ for the subgroupoid $H_a\in L_G(U_a)$. Let the
topology of $X^s$  be generated by the $M_{x,a}$ , $x\in U_a$ and
the open sets of $X$. Then its basic open sets are any set of the
form $U\cap M_{x,a}$ where $U$ is open in $X$, thus this topology
is the coarsest for which the original open sets as well as
transitive component for $H_a$ are open, and $X^s$ is
topologically the disjoint union of the transitive component for
$H_a$, each of them with its subspace topology from $X$. Since the
topology on $X^s$ is finer than that of $X$, $I\colon
X^s\rightarrow X$, the identity map, is continuous. Hence $X^s$ is
a {\it topological foliation} . The notion of topological
foliation was defined by Ehresmann \cite{Eh1}.
\begin{thm}\label{2314}
Let $s$ be a coherent local subgroupoid of the groupoid $G$ on $X$
given by an atlas $\{(U_x, H_x): x\in X\}$.
Then the transitivity components of $glob(s)$
are connected components of $X^s$.
\end{thm}
\begin{pf}
Let $H = glob(s)$. Since $s$ is coherent ($s\leq loc(H) $), for each $a\in X$,
choose an open neighbourhood $W_a\subseteq U_a$ such that $H_a|{W_a}\subseteq H\mid_{W_a}$.
If $M$ is a transitivity component of $G$, we shall show that
\[        M = \bigcup_{a\in M}(M_{a,a}\cap W_a).       \]
If $z\in M_{a,a}\cap W_a$ for some $a\in M$, then $h\in H_a(z,a)|{W_a}$ and hence $h\in G$.
Since $a\in M$ and $M$ is the transitivity component of $G$, then $z\in M$. Hence
$M$ is a union of open sets in $X^s$ and so is open.

We prove $M$ is closed in $X^s$. Let $x\in {\overline M}$, closure is relative to $X^s$,
 then $M_{x,x}\cap W_x$ meets $M$. Let $a\in (M_{x,x}\cap W_x)\cap M$. Then,
$k\in {H_x}\mid_{W_x}(a,x)\subseteq G\mid_{W_x}$. Since $a\in M, x\in M$.
Thus $M = \overline M$ and $M$ is closed in $X^s$.

Since $M$ is clopen, if it is transitivity connected, we have to show that
it is a connected component. Since $s$ is coherent on $X$ and the
topology of $X^s$ is finer than that of $X$,
it follows that $s\leq loc(H)$ in $X^s$, from which it follows that
$glob(s)\leq glob(loc(H))$ and hence $glob(s) = glob(loc(glob(s))$, i.e.,
$H = glob(s)$ is coherent on $X^s$. Since its transitivity components are
closed by Theorem \ref{2312}, they are connected.
\end{pf}

Let $X$ be a foliated manifold which is defined by
submersions $f_i \colon U_i \rightarrow  \Reno^q$, where $\{U_i : i\in X\}$
is an open cover. Let $s$ denote the associated local equivalence relation.
If it is a $p$-dimensional foliation, then by declaring $f^{-1}(r)$
open for all $i$ and all $r\in \Reno^q$, we obtain a $p$-dimensional
manifold structure on $X$. The $f^{-1}(r)$  are called the leaves of
the foliation and the identity map becomes an immersion. Our new manifold
is exactly the space $X^s$ and  $I\colon X^s\rightarrow X$
is an immersion. This is how Bourbaki \cite{Bo1} defines a foliated manifold, by
specifying $X^s$ and the immersion $I\colon X^s\rightarrow X$.

If our foliation is given by a single submersion $p\colon X\rightarrow \Reno^q$,
then $s = loc (G)$ , where $g\in G, g=(x,y)$ iff $p(x) = p(y)$. The
groupoid $G$ is
the same as the relation of being in the same
connected leaf. Thus,
by Theorem \ref{2312},  $G =glob(loc(G))$ and
$s = loc(G)$ is globally coherent.
If $U$ is a non empty open set in $X$, then $f\colon U\rightarrow \Reno^q$ is still
a submersion and by the above arguments it is easy to see
$s\mid U\leq loc(glob(s))\mid_U$, and $s = loc(G)$  is globally and
totally coherent. Thus, for our general foliation and hence by
Theorem \ref{2314}
$s$ is totally coherent.

As a result, some properties of local equivalence relations
 can be described by the results
centred around the notion of local subgroupoids. The interplay of the functors
{\em glob} and {\em loc}  says a lot about  local subgroupoids on
arbitrary topological spaces, (for more examples, see \cite{BIM}).

\chapter{Holonomy groupoid}

We quote the following historical remark from  \cite{Ao-Br}.

The concept of holonomy groupoid was  introduced  by  C.Ehresmann
and Weishu Shih in 1956 \cite{Eh2} and C.Ehresmann in 1961 \cite{Eh1},
for  a locally simple topological foliation on a topological space
$X$ (this means that  $X$ has two comparable topologies, and with
respect to the  finer  topology  on $X$,  a cover by open sets, in
each of which the  two  topologies  coincide). Such a holonomy
groupoid is considered as a topological groupoid  $H$    on $X$.
It is constructed as a groupoid of local germs of the groupoid
$H'$ of holonomy isomorphisms between the transverse spaces  $U_i$
of simple  open subsets  $U_i$   of  $X$  such that  $(U_i
,U_{i+1})$ is a  `pure  chain'.   The holonomy group at  $x\in X$
is the vertex group   $H(x)$   of   $H$.   This holonomy group is
isomorphic to the holonomy group  $H(y)$  for each   $y$   on the
same leaf of the foliation as  $x$.

Pradines \cite{Pr1} considered this holonomy groupoid  $H$,
in a wider context, with its differential structure.  He  took
the point  of view that a foliation determines an equivalence
relation $R$   by   $xRy$   if and only if  $x$  and  $y$  are on
the same leaf of the  foliation,  and  that this equivalence
relation should be regarded as a groupoid in the  standard way,
with multiplication $(x,y)(y,z)=(x,z)$  for  $(x,y), (y,z)\in  R$.
This groupoid is also written   $R$.   In  the  paracompact  case,
the locally differential structure which gives the foliation
determines a differential structure, not on  $R$  itself, but
`locally' on  $R$, that  is,  on  a  subset $W$  of  $R$
containing the diagonal  $\Delta X$ of  $X$. That is, the
foliation determines a locally topological groupoid. The full
details of this are given in \cite{Br-Mu2}.

This led Pradines to  a  definition  of  ``un  morceau
differentiable de groupoide"  $G$, for which \cite{Ma}, p.161, uses the
term ``locally differential groupoid".  Pradines' note \cite{Pr1}
asserts essentially that  such  a   $(G,W)$ determines a
differential groupoid  $Q_0 (G,W)$  and  a  homomorphism
$P:Q_0
(G,W)\to G$  such that the ``germ"  of   $W$    extends  to  a
differential structure on $G$ if and only if  $P$   is  an
isomorphism.   However  his statement of results assumes that the
base  $X=O_G$   is paracompact and that $(G,W)$   is
$\alpha$
-connected.
These assumptions seem to  be necessary  to  extend the
Globalisation Theorem 2.1 in \cite{Ao-Br} to the case of germs.

The groupoid  $Q_0 (G,W) $ is called by Pradines the {\it holonomy
groupoid} of  $(G,W)$.

A construction of the holonomy groupoid in the  differential  case
is attempted by Almeida in \cite{Al}, using properties  of  integration
of  vector fields.  However this construction has not been
published elsewhere, and of course does not extend to the
topological case.

Following Ehresmann's work,  there  has  long  been  interest  in
the holonomy group of a leaf of a smooth foliation.
For the locally differential groupoid
corresponding to a smooth  foliation, the vertex groups of  the
Ehresmann-Pradines  holonomy  groupoid  are  the holonomy groups
in the standard sense.

The holonomy groupoid  $H$   of a smooth foliation on  a  manifold
$X$ was rediscovered  (using  a  different,  but  equivalent,
description)  by Winkelnkemper \cite{Wi}, as the ``graph of the
foliation".  This  was  defined  as the set  $S$  of all triples
$(x,y,[\gamma ])$,  where  $x,y\in X$  are  on  the same leaf  $L$
of the foliation,  $\gamma$   is a continuous path on   $L$   and
$[\gamma ]$  is the equivalence class  of   $\gamma$    under  the
equivalence relation  $\sim$   which is given by:  for the two
paths $\gamma _1$, $\gamma_2$ in  $L$  starting at  $x$  and
ending at $y$, $\gamma _1 \sim \gamma _2$   if  and only if the
holonomy of  $L$ at $x$ along  $\gamma _1 \io\gamma_ 2$  is zero.
As pointed out above, these ideas are a special case of the
general construction considered here.      The way in which the
holonomy and monodromy are related in the general case is
discussed in \cite{Br-Mu1}.

Connes \cite{Co1} has considered this differential  holonomy  groupoid
$H$ of the foliation and applied to it his general theory of
integration  based on transverse measures on a measurable
groupoid. More recently, in \cite{Co2}, he has applied this and
other groupoids in the theory of non commutative $C^*$-algebras.

Pradines in \cite{Pr1} also defines what  he  calls  a  germ  of  a
locally differential  groupoid,  by  saying  two  locally
differential groupoids $(G,W)$    and   $(G,W')$   are  {\it
equivalent}  if there is  a  third  locally differential groupoid
$(G,W'')$  such that  $W''$ is  an  open  submanifold  of both
$W$   and  $W'$. The equivalence classes form the {\it germ}
of $(G, W)$. Such a germ is called a  microdifferential
groupoid. His aim is then to define the  holonomy  groupoid  as  a
functor  on  the category of such microdifferential groupoids.
One of the problems of  this theory is  that  if   $(G,W)$    and
$(G,W')$   are locally differential groupoids, then  $W\cap W'$
may no longer generate $G$.   This difficulty does  not  occur  if
the locally differential groupoids   are    $\alpha $-connected,
since in this case if  $W$ generates  $G$   and  so also does any
open subset of $W$ containing  $O_G $.  Thus there is  still  work
to be done in investigating examples of these  constructions  and
the relations between and consequences of various possible
definitions.

Three  principal  examples  of  groupoids  are  bundles   of
groups, equivalence relations, symmetry groupoids, and action
groupoids  associated with an action of a group (or more generally
groupoid) on a  set  (see  for example \cite{Br2}).   At  present,  it
seems  that  only  the  holonomy  of  an equivalence relation has
been extensively studied, namely in  the  form  of the holonomy
groups and holonomy groupoid of a smooth  foliation  (but  see
also \cite{Ro1,Ro2,Ko-Mo1,Ph,Br-Ic,BIM}).  There is presumably
considerable
potential value  in the other cases.
 The paper
\cite{BIM} gives a
new range of examples of local subgroupoids which generalise the
foliation example -- the key idea is that of {\em star path
component of the identities.}

\section{Local subgroupoids and Locally Topological \\ Groupoids}

In this section , we  obtain a holonomy groupoid for a certain local
subgroupoid by using the idea of locally topological groupoid.
Many of the idea of this section derive from in Rosenthal's
papers \cite{Ro1,Ro2}. He obtained a holonomy groupoid of a local
equivalence relation on a topological space by using the method of
Pradines \cite{Pr1} as sketched by Brown  \cite{Br1}.

The construction of the holonomy groupoid is intimately bound up
with the properties of the admissible local section of groupoid
$G$. Now, we can give the definition due to Ehresmann \cite{Eh1},
but following the notation of \cite{Ma}, with some modifications.

\begin{defn}{\rm
An {\it admissible local section } of $G$ is a function $k: U\rightarrow  G$
from an open subset $U$ of $X$ such that $k$ satisfies:

(i) \ $\alpha k(x) = x$ for all $x\in U$

(ii)\ $\beta k(U)$ is open in $X$, and

(iii)\ $\beta k$ maps $U$ homeomorphically to  $\beta  k(U)$.

The set $U$ is called the domain of $k$ and denoted by $D(k)=U$.

Let $W$ be a subset of $G$, and  suppose that $W$ has the structure of a
topological space with $X$ as a subspace.
 We say that $(\alpha, \beta, W)$
{\it has enough continuous admissible local sections }
if for each $w\in W$
there is an admissible local section $k$ of $G$ such that

(i) \ $k\alpha (w) = w,$

(ii) \ $k(U)\subseteq W$,

(iii) $k$ is continuous as a function $U\rightarrow W$.

Such a $k$ is  called a
{\it continuous admissible local section through $w$}.}
\end{defn}
If $(\alpha, \beta, W)$
 has enough continuous admissible section,
then $(\alpha, \beta, W)$
is called {\it locally sectionable}.

The holonomy groupoid will be constructed for a locally topological groupoid,
a term we now define. This definition is a modification of one due to
J.Pradines \cite{Pr1} under the name
`{\it un morceau diff\'{e}rentiable de groupoide'.}

\begin{defn}{\rm
A {\it locally topological groupoid } is a pair $(G, W)$ consisting of
a groupoid $G$ and a topological space $W$ such that

$(G_1)$  $X = Ob(G) \subseteq W  \subseteq G$

$(G_2)$  $W = W^{-1}$

$(G_3)$  the set $W_{\delta} = (W\times_{\alpha} W)\cap \delta^{-1}(W)$
is open in $W\times_{\alpha}W$ and the restriction to $W_{\delta}$ of the difference map
$\delta : G\times_{\alpha} G\rightarrow G$ , $(g, h)\mapsto gh^{-1}$,
is continuous.

$(G_4)$  the restriction to $W$ of the source and target maps $\alpha$ and $\beta$
are continuous, and the triple $(\alpha, \beta, W)$ is locally
sectionable.

$(G_5)$ $W$ generates $G$ as a groupoid.}

\end{defn}

Note that, in this definition, $G$ is a groupoid but does not need to have
a topology. In the cases considered later, $G$ will be a topological groupoid,
and $W$ is a subspace, so that condition $G_3)$ is automatic. However, $W$
will usually not be open in $G$.

Now, we will extend some definitions from local equivalence relations,
as given in Rosenthal \cite{Ro2}, to the local subgroupoids
defined in the
previous chapter.

\begin{defn}
\label{323}{\rm
Let $s$ be a local subgroupoid of a topological groupoid $G$
on $X$. An atlas $\{(U_x, H_x): x\in X, H_x\in L_G(U_x)\}$
 is called {\it weakly s-adaptable} if

(i)\ the atlas $\{(U_x, H_x): x\in X, H_x\in L_G(U_x)\}$
locally defines $s$.

(ii)\ glob(s) is the subgroupoid of $G$ generated
 by $\{H_x\}_{x\in X}$.
}
\end{defn}

The definition of $r$-adaptable atlas, (he really in fact used the
term  `family')
 which is due to Rosenthal \cite{Ro2}
for the case of equivalence relation, includes one more condition.
This is the $\alpha$-connectedness, i.e., each equivalence relation
$H_x$ has connected equivalence classes.
Then   our new the definition leads us to definition
of $\alpha$-connected locally topological groupoid \cite{Ao-Br}. That is why we do not consider
this condition.
\begin{prop}
Let $s$ be a totally coherent local subgroupoid of the topological groupoid $G$
on $X$. Then $s$ admits a weakly $s$-adaptable atlas.
\end{prop}
\begin{pf}
By assumption, we may suppose  $s$ is   defined by  an atlas
$\{(U_x,K_x):x\in
X\}$. By coherence, there is an open
cover ${\cal V}= \{V_x\}_{x\in X}$ of $X$
with  $x\in V_x\subseteq U_x$ such that $glob(s)$ is the a groupoid   $K_{\cal V}$
generated by $\{K_x\mid_{V_x}\}_{x\in X}$. By corollary
\ref{2311},
$s|{V_x}$ is globally coherent hence

\[    s|{V_x} = loc(glob(s|{V_x})).                       \]

Since $glob(s) = K_{\cal V}$ and in  some neighbourhood
of $x$, we have $K_x = H_x = glob (s|{V_x})$.
It follows that Definition \ref{323},(ii) above will
be satisfied for $(K_x, V_x)$
\end{pf}

We emphasise  that the relationship between this work and Rosenthal's is
that an equivalence relation on $X$ is simply a wide subgroupoid of the
groupoid $X\times X$, and we need $X$ to be a topological space to define
a local equivalence relation on $X$. So the appropriate content for
this work seems to be that of local subgroupoid of a topological
groupoid $G$ on $X$.

\begin{defn}
{\rm
A local subgroupoid $s$ is said to be {\it regular} if it is totally
coherent and has a weakly $s$-adaptable atlas
$\{(U_x, H_x): x\in X\}$
such that for all $x\in X$, $(\alpha_x, \beta_x, H_x)$  is
locally sectionable.
}
\end{defn}

\begin{defn}
{\rm
Let $s$ be a local subgroupoid of the topological groupoid $G$
on $X$. Then we say that $s$ is a {\it strictly regular }
local subgroupoid
if it has a regular weakly $s$-adaptable
atlas
$\{(U_x, H_x): x\in X\}$
such that, for each $g\in H_x(x, z)$ and $h\in H_y(x, y)$
, then $gh^{-1}\in H_z(y, z)$.
}
\end{defn}

We now give a key construction of a locally topological groupoid from a
strictly regular local subgroupoid.

\begin{thm}\label{327}
Let $G$ be a topological groupoid on $X$ and $s$ be a
strictly regular local
subgroupoid of $G$ on $X$ defined by atlas
${\cal U}_s=\{(U_x, H_x): x\in X\}$. Let

\[           H = glob(s)    \  \  \  \     W=W({\cal U}_s) =\bigcup_{x\in X} H_x.    \]

Then $(H, W)$ admit the structure of a locally  topological groupoid.
\end{thm}

\begin{pf}

$(G_1)$\ Because of the definition of $H$ and $W$, clearly  $X\subseteq W\subseteq H$.

$(G_2)$ In fact, $W = W^{-1}$. Let $g\in W$. Then there is an element $x\in X$ such that
$g\in H_x$. Since $H_x$ is a groupoid on $U_x$, \ $g^{-1}\in H_x$. So $W = W^{-1}$.

$(G_3)$\ We will show that $W_{\delta} = (W\times_{\alpha} W)\cap \delta^{-1}(W)$
is an open subset in $W\times_\alpha W$.
We have to show that, for a base open set
$U\times V$ in $G\times_\alpha G$,
\[             (U\times V)\cap (W\times_\alpha W)\subseteq \delta^{-1}(W) .     \]
Let $(k,l)\in (U\times V)\cap (W\times_\alpha W)$. Then  $(k,l)\in W\times_\alpha W$.
By the definition of $W$, there exist $x, y\in X$, $k\in H_x(x,z), \ l\in H_y(x, y)$.
Since $s$ is strictly regular, $kl^{-1}\in H_z(y,z)$. This shows that
$(k, l)\in \delta^{-1}(W)$. Hence $W_\delta$ is an open set in $W\times_\alpha W$.

We now prove the restriction
of $\delta$ to $W_{\delta}$ is smooth.
Since $G$ is a topological groupoid, for each $x\in X$, $H_x$ is a topological groupoid on $U_x$ and
so the difference map
\[               \delta_x : H_x\times H_x \rightarrow H_x             \]
is continuous. Because  $H_x\subseteq W, \ x\in X$, using the continuity of the inclusion
map $i_x: H_x\rightarrow U_x$, we get a continuous map

\[         i_x\times i_x : H_x\times_\alpha H_x\rightarrow W\times_{\alpha} W    \]
the restriction of $W_\delta $ is also continuous, that is,
\[         i_x\times i_x : H_x\times_\alpha H_x\rightarrow W_\delta    \]
is continuous. Then the following diagram is commutative;
$$
\begin{diagram}
              H_x\times_{\alpha_x} H_x   &\rArr        & H_x     \\
               i_x\times i_x\dArr        &             &\dArr i_x  \\
                              W_{\delta} &\rArr        &W
\end{diagram}
$$
This verifies $(G_3)$, since $H_x$ is open in $W$  and hence
$H_x\times_{\alpha} H_x$
 is open $W_\delta$.

$(G_4)$ We define source and target maps $\alpha_W$ and $\beta_W$ respectively as
follows: if $g\in W$ there exist $x\in X$ such that $g\in H_x$  and we let
\[      \alpha_W(g) = \alpha_x(g) \  \   \beta_W(g) = \beta_x(g)       \]
Clearly $\alpha_W$ and $\beta_W$ are continuous.
Since $\{(U_x, H_x): x\in X, H_x\in L_G(U_x)\}$
is a regular weakly $s$-adaptable atlas,  so
$(\alpha_x, \beta_x, H_x)$
is locally sectionable, for all $x\in X$.
Hence $(\alpha_W, \beta_W, W)$ is locally sectionable.

$(G_5)$  By definition of weakly $s$-adaptable atlas, $glob(s)$ is a subgroupoid which is generated
by $\{H_x\}_{x\in X}$, then $W$ generates $H$.

Hence $(H, W)$ is a locally topological groupoid.
\end{pf}

The following basic example  is given in Brown-Mucuk  \cite{Br-Mu2}:

Let $X$ be a foliated paracompact manifold and
let $\{U_x: x\in X\}$ be a distinguished chart
of $X$. We write $R_x$ for the equivalence
relation on $U_x$ given by $uR_xv$ if $u, v$ belong to the same path component
of $U_x$ with the leaf topology, i.e., $u$ and $v$ are in the same leaf. This
equivalence relation defines a local equivalence relation $s$ on $X$.
So we can get  an atlas $\{(U_x, R_x): x\in X\}$
which locally defines $s$ and
$glob(s) = R$. The atlas $\{(U_x, R_x): x\in X\}$
 is  weakly $s$-adaptable.
Let
\[               W = \bigcup_{x\in X}R_x                 \]
and let $W$ have its topology as a subspace of $X\times X$. Then $W\subseteq R$
but in general $W$ is not open in $R$. The triple $(\alpha, \beta, W)$
has enough continuous admissible sections. So
$(\alpha_x, \beta_x, R_x)$  has enough continuous admissible sections.
Hence $s$ is a regular local equivalence relation on $X$.
However the strictly regular  condition is proved by using the
distinguished chart and paracompactness.

\section{Holonomy groupoid}

There is a main globalisation theorem for a locally
topological groupoid. Aof-Brown in \cite{Ao-Br} stated this main theorem,
which shows how a locally topological
groupoid gives rise to its holonomy groupoid, which
is a topological
groupoid satisfying a universal property.
This theorem generalises
Th\'{e}oreme 1 of Pradines \cite{Pr1}.
Now we state  this theorem
 for  certain local subgroupoids.

\begin{thm}
Let $s$ be a local subgroupoid of a topological space $G$ on $X$,
and suppose given a strictly regular atlas
$\{(U_x, H_x): x\in X, H_x\in L_G(U_x)\}$  for $s$.
Let $(H, W)$ be  the associated locally topological groupoid.
 Then there is a topological
groupoid $Hol^s$, a morphism $\phi : Hol^s\rightarrow H $ of groupoids and an
embedding $i : W\rightarrow Hol^s$ of $W$ to an open neighbourhood of $Ob(Hol^s) = X$
such that the following condition are satisfied.

(i) $\phi $ is the identity on objects, $\phi i = id_W, \  \phi^{-1}(W)$ is
open in $Hol^s$, and the restriction $\phi_W : \phi^{-1}(W)\rightarrow W$
is continuous.

(ii) if $A$ is a topological groupoid and $\psi :A\rightarrow H$ is a morphism
of groupoids such that

(a) $\psi$ is the identity on objects,

(b) the restriction $\psi_{H_x}:\psi^{-1}(H_x)\rightarrow H_x$ of
$\psi$ is continuous and $\phi^{-1}(H_x)$
is open in $A$, the union of $\phi^{-1}(H_x)$   generates $A$

(c) the triple $(\alpha, \beta, A)$ is locally sectionable.

Then there is a unique morphism $\psi ':A\rightarrow Hol^s$ of topological groupoids such that
$\phi \psi ' = \psi$ and $\psi 'a = i\psi a $ for $a\in \psi^{-1}(W)$.
\end{thm}
The groupoid $Hol^s$ is called the holonomy groupoid $Hol^s(H, W)$ of
the local subgroupoid $s$.

We now give the construction of holonomy groupoid as in Aof-Brown  \cite{Ao-Br}.
Let $H=glob(s)$.
Let  $\Gamma(H)$ be the set of all local
admissible sections of $H$. Define a product on $\Gamma(H)$ by

\[                (tk)(x) = t(\beta k(x))k(x)    \  \   \   \   (*)         \]
for two admissible local sections $k$ and $t$. If $k$ is an
admissible local section then write $k^{-1}$ for the admissible
local section $\beta k(U)\rightarrow H$, $\beta k(x)\mapsto
(k(x))^{-1}$. With this product  $\Gamma(H)$ becomes an inverse
semigroup. Let $\Gamma^c(W)$ be the subset of  $\Gamma(H)$
consisting of admissible local sections which have values in $W$
and are continuous. Let $\Gamma^c(H, W)$ be the subsemigroup of
$\Gamma(H)$ generated by $\Gamma^c(W)$. Then $\Gamma^c(H, W)$ is
again an inverse semigroup. Intuitively, it contains information
on the iteration of local procedures. Let $J(H)$ be the sheaf of
germs of all admissible local sections of $H$. Thus the elements
of $J(H)$ are equivalence classes of pairs $(x,k)$ such that $k\in
\Gamma (H)$ , $x\in U = D(k)$  and $(x,k)$ is equivalent to
$(y,t)$ if and only if  $x = y$ and $k$ and $t$ agree on a
neighbourhood of $x$. The equivalence class of $(x,k)$ is written
$[k]_x$. The product structure on $\Gamma(H)$ induces a groupoid
structure on $J(H)$ with $X$ as the set of objects and source and
target maps $[k]_x\mapsto  x, \   [k]_x\mapsto \beta k(x)$. Let
$J^c(H, W)$ be generated as a subgroupoid of $J(H)$ by the sheaf
$J^c(W)$ of germs of element of $\Gamma^c(W)$.

Thus an element of $J^c(H,W)$ is of the form
\[            [k]_x = [k_n]_{x_n}...[k_1]_{x_1}          \]
where $k = k_n,...,k_1$ with $[k_i]\in J^c(W)$,
$x_{i+1} = \beta k_i(x_i) \
i =1,...,n$ and $x_1 = x\in U = D(k)$.

Let $\psi :J(H)\rightarrow H$ be the final map defined by $\psi
([k]_x) = k(x)$, where $k$ is an admissible local section.Then
$\psi(J^c(H,W)) = H$ because $W$ generates $H$.  Let $J_o = J^c(W)
\cap Ker\psi $. Then $J_o$ is a normal subgroupoid of $J^c(H,W)$.
The holonomy groupoid $Hol^s = Hol(H,W)$ is defined to be the
quotient groupoid $J^c(H,W)/J_o$. Let $p : J^c(H,W)\rightarrow
Hol^s$ be the quotient morphism and let $p([k]_x)$ be denoted by
$\<k\>_x$. Since $J_o\subseteq Ker\psi$ there is a surjective
morphism $\phi : Hol^s\rightarrow H$ such that $\phi p = \psi$.

The topology on the holonomy groupoid $Hol^s$ such that $Hol^s$
with this topology is
a topological groupoid is constructed as follows.
Let $k\in \Gamma^c(H,W)$ with domain $U$.
A partial function $\sigma_k:W\rightarrow Hol^s$ is defined as follows.
The domain of $\sigma_k$ is the set of $w\in W$ such that
$\beta (w)\in U$.
A continuous admissible local section  $f$ through $w$ is
chosen and the value
$\sigma_k(w)$ is defined to be
\[           \sigma_k(w) = \<k\>_{\beta (w)}\<f\>_{\alpha (w)} = \<kf\>_{\alpha w}         \]
Now we prove a Lemma which shows  that $\sigma_k (w)$ is independent of the choice of
the local section $f$.
\begin{lem} \label{332}  Let  $w\in W$, and let  $s$  and  $t$  be continuous
admissible local sections through  $w$.  Let  $x=\alpha w$.  Then
$\<s\>_x  =\< t\>_x$ in  $H$. \end{lem}

\begin{pf}  By assumption   $sx=tx=w$.   Let   $y=\beta  w$.
Without loss   of generality, we may assume that  $s$  and  $t$
have the same domain $U$   and have image contained in   $W$.
Clearly $st\io \in \Gamma^c   (W)$.   So $[st\io  ]_y \in J_0 $.
Hence in $H$ $$\< t\>_x =\< st^{-1}\>_y \< t\>_x =\<s\>_x .$$
\end{pf}
It is proven that $\sigma_k (w)$ is independent of the choice of
the local section $f$ and that these $\sigma_k$  form a set of charts.
Then the initial topology
with respect to the charts $\sigma_k$ is imposed on $Hol^s$. With this topology $H$
becomes a topological groupoid. The proof is essentially the same as
in Aof-Brown \cite{Ao-Br}.

Note that recently the structure given above have been extensively
generalised to Lie local subgroupoid and their holonomy and
monodromy Lie groupoid \cite{Br-Ic}.

\chapter{ s-sheaves }

A central area in the applications of  topology is the relation
between local and global properties. Many ideas have been developed
for this, including cohomology, sheaves and spectral sequences.
More recently, groupoids have played an increasing role.

Our attention will be focused on the interaction of sheaves and
groupoids. In this Chapter, we consider the notion of local
subgroupoid of a topological groupoid $G$, which is a global section
of a certain sheaf of subgroupoids associated to $G$. We  then
construct a category of action of a local subgroupoid.

The background to this idea comes from the notion of
local equivalence
relation and their relation with certain topos of sheaves called
\'etendues, which are categories of sheaves equipped with an action of
 an \'etale topological groupoid $G$. The aim of this chapter is
to generalise $r$-sheaf which obtained 
from strictly regular open local equivalence relation,  to the local
subgroupoids case.

To understand all structure from the beginning, now we shall give
some basic definitions  which  are due to Mac Lane and Moerdijk
\cite{Ma-Mo}.

\section{Internal Category}

Let $K$ be a category with pullbacks. Just as an ordinary small category consists of a set of objects
and a set of morphisms, so an {\it internal category} $C$ in $K$ consist
of two objects of $K$ -an  `object of objects' $C_0$ and 'object of morphisms'
$C_1$, together with four arrows of $K$ an arrow $m$ for composition, and
three arrows
$$
\begin{diagram}
C_1\times_{C_0}{C_1}  & \rArr^{m}  & C_1 & \pile{\rArr^{\alpha ,\beta}\\ \rArr \\
\lArr_{i} } & C_0
\end{diagram}
$$
for domain $\alpha$, codomain $\beta$, and identities $i$: with the first
two, we define the object $C_2$ of `compososable pairs' of morphisms as
the pullback
$$
\begin{diagram}
           C_2 = C_1\times_{C_0} C_1 & \rArr^{\pi_2} & {C_1} \\
             \dArr^{\pi_2}  &                     & \dArr_{\beta}  \\
            {C_1}           &\rArr_{\alpha }         & {C_0}.
\end{diagram}
$$

Indeed,  a generalized element  $h\colon X\rightarrow C_1\times_{C_0} C_1$
is thus just a pair of such elements $f, g \colon X\rightarrow C_1$
with $\alpha f = \beta g $, that is,  `a composable pair'. We now
require, in addition to the morphisms in above diagram, a fourth morphism
in $K$
\[   m : C_2 = C_1\times_{C_0} C_1\rightarrow C_1       \]
to represent composition of composable pairs. The axioms for an internal
category then require, besides the usual identities
$\alpha i = \beta i = 1 $ and $\alpha m = \alpha \pi_2$, $ \beta m =
\beta \pi_1$, commutativity of the following two diagrams which express
the associative law and the  unit law for composition;
$$
\begin{diagram}
         C_1\times_{C_0}C_1\times_{C_0}\times C_1 & \rArr^{1\times m}  & C_1\times_{C_0}C_1 \\
         \dArr^{m\times 1}                        &                    &\dArr_{m} \\
         C_1\times_{C_0} C_1                      & \rArr_{m}          & C_1
\end{diagram}
$$
$$
\begin{diagram}
             C_1\times_{C_0} C_0 & \rArr^{1\times i} & C_1\times_{C_0} C_1 & \rArr^{i\times 1} & C_0\times_{C_0} C_1  \\
                                 & \SE_{\pi_1}       & \dArr_{m}         &  \SW_{\pi_2}        &   \\
                                 &                   & C_1               &                     &
\end{diagram}
$$
These conditions constitute a `diagrammatic' form of the standard
definition
of a category. If $C$ and $D$ are two such internal
categories, then an
internal functor $F\colon C\rightarrow D$ is
defined to be a pair of
morphisms $F_0 : C_0\rightarrow D_0$ and
$F_1 :C_1\rightarrow D_1$ in $K$
making the obvious four squares ( with $i$, $\alpha$, $\beta$ , $m$)    commute.

With the evident composition of such functors
we have a category $Cat (K)$ with the internal categories in $K$ as objects
and internal functors as morphism.

In ordinary category theory, the functors $F :C\rightarrow D$
between two small categories play a role quite different from
functors from $C$ into the ambient category - the category of
sets. A functor of the latter sort consists as usual of an `
object function' $C_0\rightarrow Sets$ and an arrow 'functions'
$C_1 \rightarrow Functions$, suitably related. In this
description, we now wish to replace $Sets$ by any category $K$
with pullbacks, and $C$ by the internal category (again called
$C$) in $K$. In order to get a suitable `internal' description of
such functors to the universe $K$, we first reformulate the usual
case where the universe is $Sets$. There an  object function $F_0
:C_0\rightarrow Sets$ can be viewed as a $C_0 -indexed$ family of
sets, one for each $x\in C_0$. Just as in the treatment of indexed
sets, this $C_0-indexed$ family can be replaced by a single object
over $C_0$.
\[     p : F\rightarrow C_0     \]
where $F = \bigcup_{x\in C_0} F_0(x)$ is the disjoint sum of all
the sets $F_0(x)$, and $p$ is the obvious projection. Each set
$F_0x$ can then be recovered ( up to isomorphism ) from $p$ as the
fiber $P^{-1}(x)$. Similarly, for the arrow function, each arrow
$f:x\rightarrow y$ in $C$ given a map $F_0x\rightarrow F_0y$ of
sets, written for $a\in F_0x$ as $a\mapsto f.a$. All these maps,
one for each $f\in C_1$, can be described in terms  of $p
:F\rightarrow C_0$ as one single map specifying the action of any
$f$ on any a as
\[    \phi : C_1\times_{C_0} F\rightarrow F, \ \  \  \phi(f, a) = f.a   \]
where $C_1\times_{C_0} F\rightarrow F$ is pullback of $p$ along
$\alpha : C_1\rightarrow C_0$.

By writing down the appropriate diagrams, the preceding description of
a functor to $Sets$ can be easily generalized to the case of an interval
category $C$ in a category $K$ with pullbacks. A {\it C-object} in $K$
(also called an `internal diagram' on $C$) is an object
$p :F \rightarrow C_0$ over $C_0$ equipped with an {\it action }
\[      \phi : C_1\times_{C_0} F\rightarrow F     \]
of $C$ of $F$, where for this pullback $\alpha : C_1\rightarrow C_0$ is used
to take $C_1$ an object over $C_0$. Here the following diagrams are required
to commute:
$$
\begin{diagram}
           C_1\times_{C_0} F & \rArr^{\phi} & F \\
             \dArr^{\pi_1}  &                     & \dArr_{p}  \\
            {C_1}           &\rArr_{\beta }         & {C_0}.
\end{diagram}
\begin{diagram}
             C_0\times F &  \rArr^{i\times 1} & C_1\times_{C_0} F  \\
                         & \SE_{\pi_2}        & \dArr_{\phi}          \\
                         &                    &  F
\end{diagram}
$$
$$
\begin{diagram}
           C_1\times_{C_0} C_1\times_{C_0} F &\rArr^{1\times \phi}  & C_1\times_{C_0} F \\
           \dArr^{m\times 1}                 &                       &\dArr^{\phi}   \\
           C_1\times_{C_0} F                 &\rArr^{\phi}           &  F
\end{diagram}
$$ (The second and third express the unit and associativity laws
for the action.)

If  $F = ( F, p, \phi_1)$ and $G = (G, q, \phi_2)$ are two such
$C$-object
in $K$, {\it a morphism} of $C$-objects from $F$ to $G$ is simply a
morphism
$\psi : F\rightarrow G$ in $K$ which preserves the
structure involved.
In term of diagram, it means that
$$
\begin{diagram}
            F & \rArr^{\ \ \ \eta} &    &        &     G \\
              & \SE_{p}      &    &\SW_{q} & \\
              &              &  C_0 &        &
\end{diagram}
\begin{diagram}
\  \  \  \             C_1\times_{C_0} C_1\times_{C_0} F &\rArr^{1\times \phi}  & C_1\times_{C_0} F \\
\  \  \  \              \dArr^{m\times 1}                 &                       &\dArr{\phi}   \\
\  \   \  \           C_1\times_{C_0} F                 &\rArr_{\phi}           &  F
\end{diagram}
$$
are required to commute.

%---------------------------------------------------------
\section{G-sheaves}

Let
$$
\begin{diagram}
{\bf C } =  & C_2 &\pile{\rArr^{\pi_1} \\ \rArr^{\pi_2} \\ \rArr^{m}}
& C_1&\pile{\rArr^{\alpha } \\ \rArr^{\beta }} & C_0        \     \   \  \  \  (\star)
\end{diagram}
$$
be a topological internal  category, where $C_0$ is the space of  objects, $C_1$  is
the space of morphism and $C_2$ is  the  space  of   composable
pairs  of
morphisms, $m$ represents composition, $\alpha $ and $\beta $  are
the domain and codomain
maps, respectively,  which have a common section $i \colon C_0\rightarrow C_1$.
\begin{defn}{\rm

 Let {\bf C} be a topological category. A {\bf C}-Sheaf  is  a  sheaf
$  p\colon {\cal F}\rightarrow X$ together with a diagram
$$
\begin{diagram}
               { C_1\times {\cal F}} &  \rArr^{\phi } &  {\cal F}  \\
                 \dArr^{\pi_1}        &                &  \dArr_{p}  \  \  \  \ (**)   \\
                  C_1         &  \rArr_{\beta }   &C_0
\end{diagram}
$$
which is commutative and furthermore:}
\end{defn}
$$
\begin{diagram}
\   \    \   C_0\times {\cal F} & \rArr^{i \times I} & C_1\times {\cal F} \\
 \      \dArr^{\pi_2} & \SW_{\phi}              &  \\
 \            {\cal F}     &                         &
\end{diagram}
\begin{diagram}
\  \ \  \  C_1\times C_1\times {\cal F} & \rArr^{1\times \phi} & C_1\times {\cal F} \\
    \dArr^{m\times I}&                             &\dArr_{\phi } \\
    C_1\times {\cal F}  &    \rArr                           &{\cal F}
\end{diagram}
$$

The first and second express respectively  the unit and associativity laws for
the map $\phi$, i.e., $\phi $ is an action of ${\cal F}$  on $X$ \cite{Ma-Mo}.

A {\bf C} -Sheaf morphism is a sheaf map which preserves the actions.
In this way, we obtain the category $Sh(X: {\bf C})$ of {\bf C} -sheaves.

Now, let $G$ be a  topological groupoid  on $X$. If
we take  $G =C_1$ and $G\times G = C_2$ in $(\star)$
we obtain a ${\bf G}-
topological$  category . In this category, $\alpha$ and $\beta $
are source  and  target  map  , $m$
represents composition, $\pi_1$  and $\pi_2$  are the canonical  projections
and $\epsilon $ is the object map.
$$
\begin{diagram}
{\bf C } =  & G\times_X G &\pile{\rArr^{\pi_1} \\ \rArr^{\pi_2} \\ \rArr^{m}}
& G&\pile{\rArr^{\alpha } \\ \rArr^{\beta }} & X       \     \   \  \  \
\end{diagram}
$$

A sheaf $p\colon {\cal F}\rightarrow X$ is called
${\bf G}$-sheaf  if  it  satisfies (**),
unit and associativity laws   above, i.e.,
let us given following diagram as in (**),
$$
\begin{diagram}
             G\times_X {\cal F} & \rArr^{\phi} &{\cal F} \\
             \dArr^{\pi_1}       &              &\dArr_{p} \\
             G          &\rArr^{\beta}   & X
\end{diagram}
$$
where  $G \times_X {\cal F}$ is a pullback on $\beta $. Then, it
satisfies

(i) $p(\phi (g)(e_x)) = y$, \ \  for $g\in G(x,y)$

(ii) $\phi ((g_x)(e_x)) = e_x $, \ \ for $e_x\in {\cal F}_x$

(ii)  $\phi (h)(\phi (g)(e_x)) = \phi (k)(e_x)$, \ \
where $g_x\in G(x,x) , g\in G(x,y),  h\in G(y,z)$ and, $k=gh\in G(x,z)$.

So each  element $g$ in $G(x,y)$ defines a morphism
\[       g_{\sharp }\colon {\cal F}_x \rightarrow {\cal F}_y                        \]
of stalks such that   $I_{\sharp }=I$ and  $(hg)_{\sharp } = g_{\sharp } h_{\sharp }$.
Thus an action of $G$ on  ${\cal F}$ defines a  functor
\[    {\cal F}  \colon G \rightarrow \ \{ \mbox{stalks of sheaves }. \}   \]
The map $\phi $ is called  {\it a transport} along $G$ in ${\cal F}$.

We give as an example of $G$-sheaf by taking an
equivalence relation $R$ on $X$. Then
$$
\begin{diagram}
{\bf R } =  & R\times R &\pile{\rArr^{\pi_1} \\ \rArr^{\pi_2} \\ \rArr^{m}}
& R &\pile{\rArr^{\pi_1 } \\ \rArr^{\pi_2 }} & X        \     \   \  \  \  (\star)
\end{diagram}
$$
is an internal topological category. Since $R$ is a topological groupoid,
it can be shown that an  $R$-sheaf can be described as a sheaf
$ p :{\cal F}\rightarrow X$ with an equivalence relation $S$ on {\cal F}
such that if $(x_1, x_2)\in R$ and $e_1\in p^{-1}(x_1)$, there is a unique
$e_2\in p^{-1}(x_1)$ with $(e_1, e_2)\in S$.

\begin{prop}
The category of {\bf G}-sheaves is a Grothendieck topos.
\end{prop}
\begin{pf}
See  Moerdijk \cite{Mo1}.
\end{pf}
This topos is denoted by $BG$, and called the {\it classifying
topos}  of $G$. The definition of the topos $BG$ also makes sense
if $G$ is just a continuous category (a category object in
Locales) rather than a groupoid. We can find a lot of examples of
the topos $BG$ in \cite{Mo1}.

%-----------------------------------------------------
\section{r-sheaf}

Let $ p :{\cal F}\rightarrow X$  be a sheaf over $X$, and let $U$ be open in
$X$. Let $Q(U, {\cal F})$ consist of pairs $(R_U, S_U)$,  where $R_U$ and
$S_U$ are equivalence relation on $U$, ${\cal F}\mid_U$, respectively,
such that $p :{\cal F}\mid_U\rightarrow U$ is compatible with $R_U$ and
$S_U$ i.e.
\[     (e_1, e_2)\in S_U  \ \  \mbox{implies}
\ \ \ (p(e_1), p(e_2))\in R_U  \]
and we have the following pullback
$$
\begin{diagram}
           {\cal F}\mid_U  & \rArr      & ({\cal F}\mid_U) /{S_U}  \\
           \dArr^{p}       &            & \dArr_{q} \\
           U               & \rArr      &  U/{R_U}
\end{diagram}
$$
with $q$ a local homeomorphism. This implies that if
$p(e) = x$ and
$x'\in [x]$, the equivalence class of
$x$ with respect to $R_U$,
there is a unique
\[     e'\in p^{-1}(x) \ \ \ \  \mbox{with}  \ \ \ (e', e) \in S_U .    \]

We can give the following theorem.

\begin{thm}
Let {\cal F} be a sheaf on $X$.
Let  ${\cal O}(X)$ denote the open sets of
$X$. Then
\[      Q(-, {\cal F}) : {\cal O}(X)^{op}\rightarrow Sets      \]
is a presheaf.
\end{thm}
\begin{pf}
See  Rosenthal \cite{Ro2}.
\end{pf}

Let $Q_{{\cal F}}$ denote the associated sheaf. We have a forgetful
functor of sheaves $Q_{{\cal F}}\rightarrow {\cal E}$, the sheaf of local
equivalence relation on $X$.

Let $r$ be a local equivalence relation on $X$.
\begin{defn}{\rm
An $r$-structure on a sheaf ${\cal F}$ is a local equivalence relation
$t$ on ${\cal F}$ such that $(r, t)$ is a global section of
$Q_{\cal F}$, i.e., $p(t)=r$.}
\end{defn}
\begin{defn}{\rm
An $r$-sheaf on $X$ is a pair $({\cal F}, t)$, where ${\cal F}$ is a sheaf
on $X$ and $t$ is an $r$-structure on ${\cal F}$.}
\end{defn}

If $({\cal F}_1, t_1)$ and $({\cal F}_2, t_2)$ are $r$-sheaves,  an  $r$-sheaf
morphism is a sheaf map  ${\cal F}_1\rightarrow   {\cal F}_2$, which locally
preserves the $r$-structures. Thus we have a category $Sh(X; r)$ of
$r$-sheaves.

It is clear that because of $Sh(X: R) \cong Sh(X: R/X)$
\cite{Ro2}, every $R$-sheaf can be viewed as an $r$-sheaf. It is
shown that every $r$-sheaf can be made into $R$-sheaf. If $p
:{\cal F}\rightarrow X$ is an $r$-sheaf with an $r$-structure $t$,
we must produce an $R$-action $\phi : R\times {\cal F}\rightarrow
{\cal F}$ making the following diagram commute: $$
\begin{diagram}
             R\times {\cal F}  &  \rArr^{\phi} & {\cal F}  \\
             \dArr^{\pi_1}             &               &\dArr_{p}  \\
             R                 & \rArr_{\beta}          &X.
\end{diagram}
$$

\section{s-sheaf}
In this section, we shall define $s$-sheaf for a local subgroupoid
$s: X\to {\cal L}_G$ given by an atlas ${\cal U}_s=\{(U_i, H_i):i\in
I, H_i\in L_G(U_i)\}$.  Let ${\cal F}$ be a sheaf on $X$.
For any open subset $U\subseteq X$,  we consider the  set $I_{\cal F}(U)$
consisting of  pairs
$(H_i,\phi_i)$ , where  $H_i\in L_G(U_i)$  and $\phi_i $ is  a
 transport
by $H_i$ on ${\cal F}\mid_{U_i} = p^{-1}(U_i)$.
 For $U_j\subseteq U_i$, there is
a
restriction map
\[    I_{\cal F}(U_i)\longrightarrow I_{\cal F}(U_j),                \]
and this give the presheaf $I_{\cal F}$. Furthermore, there  is
a forgetful functor  $I_{\cal F}\rightarrow  L_G$ given by
$(H_i,\phi_i)\mapsto H_i$.
Let ${\cal I}_{\cal F}$ be the associated sheaf, so the
forgetful functor $I_{\cal F}\rightarrow L_G$
induces a map of sheaves
${\cal I}_{\cal F}\rightarrow {\cal L}_G$.
Fix a local subgroupoid
$s$ of $G$ on $X$.

\begin{defn}{\rm
   An {\it $s$-transport} on the sheaf ${\cal F}$ is a global  section  $t$
of ${\cal I}_{\cal F}$  such that $p(t) = s$.
An {\it $s$-sheaf}  on $X$ is a sheaf on $X$
together with  an $s$-transport.}
\end{defn}

The notation of the transport preserving map between two $s$-sheaves on $X$ can
be defined as follows;
\begin{defn}{\rm
Let $G$ be a topological groupoid on $X$ and let ${\cal F}_1$, ${\cal F}_2$ be $s$-sheaves
with transports $\phi_1$ and $\phi_2$, respectively. An $s$-sheaf morphism from
${\cal F}_1$ to ${\cal F}_2$ is a sheaf map $\eta \colon {\cal F}_1\rightarrow {\cal F}_2$ such that
the following diagram is commutative:
}
\end{defn}
$$
\begin{diagram}
             G\times {{\cal F}}_1 & \rArr^{I\times \eta } & G\times {{\cal F}}_2 \\
             \dArr^{\phi_1}  &                     & \dArr_{\phi_2}  \\
            {{\cal F}}_1           &\rArr_{\eta }         & {{\cal F}}_2
\end{diagram}
$$

Let $ Sh(X;s) $ denote the category of $s$-sheaves and
$s$-sheaf morphisms.
There exists a faithful functor from $Sh(X;s)$ to $Sh(X)$. Because every $s$-sheaf
on $X$ is a sheaf on $X$.

From the definitions, it also follows that
the property of being an $s$-sheaves  is locally property, i.e., if the
base space $X$ is covered by open sets $U$ such that the restriction
of the sheaf ${\cal F}\to X$  to each $U$ is an $s|U$-sheaf, then
${\cal F}$ is an $s$-sheaf.
\begin{cor}
Let $s$ be a locally transitive local subgroupoid of $G$ on $X$. Then
any sheaf ${\cal F}$ has at most one $s$-transport.
\end{cor}
\begin{pf}
The proof is similar to Theorem 2.2 given in \cite{Ko-Mo1}.
\end{pf}

Note that it has been shown that a local equivalence relation gives rise 
to an \'etendue, which is a particular kind of topos, and that Kock and Moerdijk 
\cite{Ko-Mo2}  have shown that every \'etendue arises from a local equivalence relation.
This suggests the problem of characterising the topos of $s$-sheaves arising from 
a local subgroupoid $s$.

\end{document}